\documentclass{amsart}

\usepackage{amsmath} 
\usepackage{amssymb}

\newtheorem{theorem}{Theorem}[section] 
\newtheorem{claim}{Claim}[theorem]
\newtheorem{lemma}[theorem]{Lemma} 
 
\newtheorem{observation}[theorem]{Observation}

\theoremstyle{definition}
\newtheorem{definition}[theorem]{Definition}

\newtheorem{problem}[theorem]{Problem}

\theoremstyle{remark}
\newtheorem{remark}[theorem]{Remark}

\numberwithin{equation}{section}
\newcommand{\forces}{\Vdash}

\newcommand{\bV}{{\bf V}} 
\newcommand{\lesdot}{\mathrel{\mathord{<}\!\!\raise 
0.8 pt\hbox{$\scriptstyle\circ$}}}

\newcount\skewfactor
\def\mathunderaccent#1#2 {\let\theaccent#1\skewfactor#2
\mathpalette\putaccentunder}
\def\putaccentunder#1#2{\oalign{$#1#2$\crcr\hidewidth
\vbox to.2ex{\hbox{$#1\skew\skewfactor\theaccent{}$}\vss}\hidewidth}}
\def\name{\mathunderaccent\tilde-3 }

\newcommand{\conc}{{}^\frown\!}
\newcommand{\lh}{{\rm lh}\/}
\newcommand{\rest}{{\restriction}}
\newcommand{\Dom}{{\rm Dom}} 
\newcommand{\Rng}{{\rm Rng}}
 
\newcommand{\proj}{{\rm proj}}
\newcommand{\suc}{{\rm succ}}

\newcommand{\rk}{{\rm rk}}
\newcommand{\vtl}{\vartriangleleft}
\newcommand{\baire}{{}^{\omega}\omega}
 
\newcommand{\bbD}{{\mathbb D}}

\newcommand{\cH}{{\mathcal H}}

\newcommand{\cF}{{\mathcal F}}

\newcommand{\bbP}{{\mathbb P}}

\newcommand{\bbQ}{{\mathbb Q}}

\newcommand{\mbR}{{\mathbb R}}

\newcommand{\bbS}{{\mathbb S}}
\newcommand{\cT}{{\mathcal T}}

\newcommand{\Gsg}{{\Game^\oplus_n}}
\newcommand{\wGs}{{\Game^\odot_n}}
\newcommand{\tGs}{{\Game^\ominus_n}}
\newcommand{\jeszcze}{{\Game^{\bar{\st}}_{\odot,_n}}}
\newcommand{\jeksi}{{\Game^{\name{\bar{\st}}^\xi}_{\odot,_n}}} 

\newcommand{\st}{{\bf st}} 

\newcommand{\vare}{\varepsilon}
\newcommand{\proc}{{\rm proc}}
\newcommand{\temp}{{\rm tmp}}

\begin{document}

\title{$n$--localization property}

\author{Andrzej Ros{\l}anowski}

\address{Department of Mathematics\\
 University of Nebraska at Omaha\\
 Omaha, NE 68182-0243, USA}

\email{roslanow@member.ams.org}

\urladdr{http://www.unomaha.edu/logic}

\thanks{The author would like to thank his wife, Ma{\l}gorzata
Jankowiak--Ros{\l}anowska for partial support of this research. He also
acknowledges partial support from the United States-Israel Binational
Science Foundation (Grant no. 2002323).} 

\date{July 2005}


\maketitle

\section{Introduction}
The present paper is concerned with the $n$--localization property and its
preservation in countable support (CS) iterations. This property was first
introduced in Newelski and Ros{\l}anowski \cite[p. 826]{NeRo93}. 
\begin{definition}
\label{nfor}
Let $n$ be an integer greater than 1. 
\begin{enumerate}
\item A tree $T$ is {\em an $n$--ary tree\/} provided that $(\forall s\in 
T)(|\suc_T(s)|\leq n)$. 
\item A forcing notion $\bbP$ has the {\em $n$--localization property\/} if
\[\forces_{\bbP}\mbox{`` }\big(\forall f\in\baire\big)\big(\exists T\in\bV
\big)\big(T\mbox{ is an $n$--ary tree and }f\in [T]\big)\mbox{ ''.}\]
\end{enumerate}
\end{definition}
In \cite[Theorem 2.3]{NeRo93} we showed that countable support products of
the $n$--Sacks forcing notion $\bbD_n$ (see Definition \ref{forcings}(1)
here) have the $n$--localization property. That theorem was used to obtain
some consistency results concerning cardinal characteristics of the ideal
determined by unsymmetric games. Soon after this, the uniform $n$--Sacks
forcing notion  $\bbQ_n$ (see Definition \ref{forcings}(2)) was introduced
in  \cite[\S 4]{Ro94} and applied in the proof of \cite[Theorem
  5.13]{Ro94}. The crucial property of $\bbQ_n$ which was used there is that
the CS iterations of $\bbQ_n$ have the $n$--localization property, but in
\cite{Ro94} we only stated that the proof is similar to that of
\cite[Theorem 2.3]{NeRo93}.  

One of the difficulties with the $n$--localization property was that there
was no ``preservation theorem'' for it. Geschke and Quickert \cite{GeQu0x}
give full and  detailed proofs of the 2--localization property for both CS 
products and CS iterations of the Sacks forcing $\bbD_2$ (and those proofs
can be easily rewritten for $n$--localization property and $\bbD_n$). And
the same proof can be repeated for $\bbQ_n$, but a more general theorem has
been missing. 

Recently, the $n$--localization property, the $\sigma$--ideal generated by
$n$--ary trees and $n$--Sacks forcing notion $\bbD_n$ have been found
applicable to some questions concerning convexity numbers of closed subsets
of $\mbR^n$, see Geschke, Kojman, Kubi\'s and Schipperus \cite{GKKS02}, 
Geschke and Kojman \cite{GeKo02} and most recently Geschke \cite{Ge05}. The
latter paper is {\em raison d'\^{e}tre} for this note --- when I read
\cite{Ge05} to write a review for {\em Mathematical Reviews\/} I wanted to
check as many technical details as I could. In \cite[\S 2]{Ge05} an
interesting forcing notion\footnote{we call it {\em the Geschke forcing\/}
here, see Definition \ref{forcings}(4)}  $\bbP_G$ was introduced and a proof
was given that it has the $n$--localization property. However, the proof
that the CS iteration of this forcing has the $n$--localization property was
left to the reader as ``similar to that for Sacks''. At first I was not sure 
about technical details of that proof, so I decided to look at $\bbP_G$ and
$\bbD_n$ together. Soon I have become convinced that a unifying theorem is
needed and this note presents a result which has such character. 

It was stated in \cite[Theorem 2.3]{NeRo93} that the same proof as for
$\bbD_n$ works also for CS iterations and products of the $n$--Silver
forcing notions $\bbS_n$ (see Definition \ref{forcings}(3)). Maybe some old
wisdom got lost, but it does not look like that {\em the same arguments work
for the $n$--Silver forcing $\bbS_n$}. As a matter of fact, we believe that
it is an open question if $\bbS_n$ and its CS iterations have the
$n$--localization property. Also, motivated by \cite[Questions 3.3,
3.4]{CRSW93} we asked if the iteration of two 2--Silver forcing notions
may add a 4--Silver real, but because of the claim in \cite[Corollary
2.4]{NeRo93} we did not state the question explicitly in the final version
of \cite{CRSW93}. In the light of what we said above, it is only proper to
pose this problem again. 

\begin{problem}
\begin{enumerate}
\item Can a finite iteration of 2--Silver forcings $\bbS_2$ add a generic
  real for the 4--Silver forcing notion $\bbS_4$?
\item Does the $n$--Silver forcing $\bbS_n$ have the $n$--localization
  property? The same about CS iterations of $n$--Silver forcings.
\end{enumerate}
{\em The author offers ``all you can drink in 3 days'' coffee/espresso in a
  place similar to {\em Caffeine Dreams\/} in Omaha for full solution to
  this problem. Partial solutions may be eligible for partial awards.}   
\end{problem}
 
It may occur that the answer to the above problem is hidden in Shelah and
Stepr\={a}ns \cite{ShSr:665}. Let us note that Remark \ref{finrem} suggests
that if we can show that finite iterations of $\bbS_n$ have the
$n$--localization property, then we will be able to handle all CS
iterations.  

The following general question remains still unsolved. 

\begin{problem}
Do CS iterations of proper forcing notions with $n$--localization property
have $n$--localization property? What if we restrict ourselves to (s)nep
forcing notions (see Shelah \cite{Sh:630}) or even Suslin$^+$ (see Goldstern
\cite{Go} or Kellner \cite{Kr0x}, \cite{Krphd})? 
\end{problem}
\medskip

\noindent {\bf Content of the paper:}\quad In the first section we introduce
several properties related to the $n$--localization property. The strongest
one, $\oplus_n$--property, does imply the $n$--localization. However not all
forcing notions around have the $\oplus_n$--property so this is why we have
weaker relatives. We also remind definitions of the forcing notions that we
are interested in and the basic facts on trees of conditions. 

The following section shows that CS iterations of forcing notions with
the $\oplus_n$--property have the $n$--localization (Theorem
\ref{iter}). Since we do not know if $\bbQ_n$ has the $\oplus_n$--property,
in the third section we somewhat weaken that property to cover more forcing
notions. From the point of view of applications Theorems \ref{seciter},
\ref{getnloc} are strongest and they include the result of the second
section. Still, we think that the proof of \ref{iter} is somewhat easier and
it is a good preparation for Section 3. 
\medskip

One should note that the proofs of our iteration theorems are very ``not
proper'' in their form. We work with trees of conditions which were used in
pre-proper era and our arguments resemble those of Ros{\l}anowski and Shelah
\cite[\S 2]{RoSh:860} and to some extend also \cite[\S A.2]{RoSh:777}.  
\medskip

\noindent {\bf Notation:}\quad Our notation is rather standard and
compatible with that of classical textbooks (like Jech \cite{J}). In forcing
we keep the older convention that {\em a stronger condition is the larger
  one}.  

\begin{enumerate}
\item $n$ is our fixed integer, $n\geq 2$. Ordinal numbers will be denoted
be the lower case initial letters of the Greek alphabet ($\alpha,\beta,
\gamma,\delta\ldots$) with possible sub- and superscripts. Natural number
will be labeled by $i,j,k,\ell,m$ (also upper cases).  

By $\chi$ we will denote a {\em sufficiently large\/} regular cardinal;
$\cH(\chi)$ is the family of all sets hereditarily of size less than
$\chi$. Moreover, we fix a well ordering $<^*_\chi$ of $\cH(\chi)$.  

\item For two sequences $\eta,\nu$ we write $\nu\vartriangleleft\eta$
whenever $\nu$ is a proper initial segment of $\eta$, and $\nu
\trianglelefteq\eta$ when either $\nu\vartriangleleft\eta$ or $\nu=\eta$. 
The length of a sequence $\eta$ is denoted by $\lh(\eta)$.

\item A {\em tree} is a family of finite sequences closed under initial
segments. For a tree $T$ and $\eta\in T$ we define {\em the successors of
$\eta$ in $T$} and {\em maximal points of $T$} by:
\[\begin{array}{rcl}
\suc_T(\eta)&=&\{\nu\in T: \eta\vartriangleleft\nu\ \&\ \neg(\exists\rho\in
T)(\eta\vartriangleleft\rho\vartriangleleft\nu)\},\\
\max(T)&=&\{\nu\in T:\mbox{ there is no }\rho\in T\mbox{ such that }
\nu\vartriangleleft\rho\}.
  \end{array}\]   
For a tree $T$ the family of all $\omega$--branches through $T$ is
denoted by $[T]$. 

\item We will consider some games of two players. One player will be
called {\em Generic\/}, and we will refer to this player as
``she''. Her opponent will be called {\em Antigeneric\/} and will be
referred to as ``he''. 

\item For a forcing notion $\bbP$, $\Gamma_\bbP$ stands for the canonical
$\bbP$--name for the generic filter in $\bbP$. With this one exception, all
$\bbP$--names for objects in the extension via $\bbP$ will be denoted with a
tilde below (e.g., $\name{\tau}$, $\name{X}$). The weakest element of $\bbP$
will be denoted by $\emptyset_\bbP$ (and we will always assume that there is
one, and that there is no other condition equivalent to it). We will also
assume that all forcing notions under considerations are atomless. 

By ``CS iterations'' we mean iterations in which domains of
conditions are countable. However, we will pretend that conditions in
a CS iteration $\bar{\bbQ}=\langle\bbP_\zeta,\name{\bbQ}_\zeta:\zeta<
\gamma\rangle$ are total functions on $\gamma$ and for $p\in\lim(
\bar{\bbQ})$ and $\alpha<\gamma$ we have $\forces_{\bbP_\alpha} p(\alpha)\in
\name{\bbQ}_\alpha$, and if $\alpha\in\gamma\setminus\Dom(p)$ then
$p(\alpha)=\name{\emptyset}_{\name{\bbQ}_\alpha}$.  
\end{enumerate}

\section{Tools}
In this section we introduce the main concepts and properties ans we show
how they are related to various forcing notions. We also introduce the main
tool for our forcing arguments: trees of conditions. 

\begin{definition}
\label{da2}
Let $\bbP$ be a forcing notion.
\begin{enumerate}
\item For a condition $p\in\bbP$ we define a game $\Gsg(p,\bbP)$ of two
  players, {\em Generic\/} and {\em Antigeneric\/}. A play of $\Gsg(p,\bbP)$
  lasts $\omega$ moves and during it the players construct a sequence
  $\langle (s_i,\bar{p}^i,\bar{q}^i):i<\omega\rangle$ as follows. At a stage 
  $i<\omega$ of the play, first Generic chooses a finite $n$--ary tree $s_i$
  and a system $\bar{p}^i=\langle p^i_\eta:\eta\in\max(s_i)\rangle$ such
  that:  
\begin{enumerate}
\item[$(\alpha)$] $|\max(s_0)|\leq n$ and if $i=j+1$ then $s_j$ is a subtree
of $s_i$ such that  
\[\big(\forall \eta\in\max(s_i)\big)\big(\exists \ell<\lh(\eta)\big)\big(
\eta\rest \ell\in \max(s_j)\big),\] 
and 
\[\big(\forall\nu\in\max(s_j)\big)\big(0<\big|\big\{\eta\in\max(s_i):\nu
\vtl\eta\big\}\big|\leq n\big),\] 
\item[$(\beta)$]  $p^i_\eta\in\bbP$ for all $\eta\in \max(s_i)$, 
\item[$(\gamma)$] if $j<i$, $\nu\in \max(s_j)$ and $\nu\vtl\eta\in
  \max(s_i)$, then $q^j_\nu\leq p^i_\eta$ and $p\leq p^i_\eta$. 
\end{enumerate}
Then Antigeneric answers choosing a system $\bar{q}^i=\langle q^i_\eta:
\eta\in\max(s_i)\rangle$ of conditions in $\bbP$ such that $p^i_\eta\leq
q^i_\eta$ for each $\eta\in \max(s_i)$.  

Finally, Generic wins the play $\langle (s_i,\bar{p}^i,\bar{q}^i):i<\omega 
\rangle$ if and only if 
\begin{enumerate}
\item[$(\circledast)$] there is a condition $q\geq p$ such that for every
  $i<\omega$ the family $\{q^i_\eta:\eta\in\max(s_i)\}$ is predense above
  $q$.  
\end{enumerate}
\item Let $p\in\bbP$. We define a game $\wGs(p,\bbP)$ of two players, {\em
  Generic\/} and {\em Antigeneric\/}. A play of $\wGs(p,\bbP)$ lasts
  $\omega$ moves and during it the players construct a sequence $\langle
  (s_i,\bar{p}^i,\bar{q}^i):i<\omega\rangle$ as follows. At a stage  
  $i<\omega$ of the play, first Generic chooses a finite $n$--ary tree $s_i$
  such that the demand $(\alpha)$ of (1) above holds. Next 
\begin{enumerate}
\item[$(\odot)$] Antigeneric picks an enumeration $\langle\eta^i_\ell:\ell
  <k_i\rangle$ of $\max(s_i)$ (so $k_i<\omega$) 
\end{enumerate}
and then the two players play a subgame of length $k_i$ alternatively
choosing successive terms of a sequence $\langle p^i_{ \eta^i_\ell},
q^i_{\eta^i_\ell}:\ell<k_i\rangle$. At a stage $\ell<k_i$ of the subgame,
first Generic picks a condition $p^i_{\eta^i_\ell}\in\bbP$ such that    
\begin{enumerate}
\item[$(\gamma)^i_\ell$] if $j<i$, $\nu\in \max(s_j)$ and $\nu\vtl
\eta^i_\ell$, then $q^j_\nu\leq p^i_{\eta^i_\ell}$ and $p\leq
p^i_{\eta^i_\ell}$,    
\end{enumerate}
and then Antigeneric answers with a condition $q^i_{\eta^i_\ell}$ stronger
than $p^i_{\eta^i_\ell}$. 

The winning criterion for the game $\wGs$ is the same as the one for
$\Gsg$ (i.e., $(\circledast)$).

\item A game $\tGs(p,\bbP)$ for $p\in\bbP$ is defined like $\wGs(p,\bbP)$
  above, but $(\odot)$ is replaced by 
\begin{enumerate}
\item[$(\ominus)$] Generic picks an enumeration $\langle\eta^i_\ell:\ell<k_i
  \rangle$ of $\max(s_i)$.
\end{enumerate}
\item We say that $\bbP$ has the {\em $\oplus_n$--property\/} whenever
Generic has a winning strategy in the game $\Gsg(p,\bbP)$ for any
$p\in\bbP$. In a similar manner we define when $\bbP$ has the {\em
$\odot_n$--property\/} ({\em $\ominus_n$--property}, respectively)
replacing the game $\Gsg$ by $\wGs$ ($\tGs$, respectively).   
\end{enumerate}
\end{definition}

\begin{definition}
\label{nice}
Let $\bbP$ be a forcing notion.
\begin{enumerate}
\item Assume that $K\subseteq\omega$ is infinite, $p\in\bbP$. A strategy
$\st$ for Generic in $\wGs(p,\bbP)$ is said to be {\em nice for $K$}
(or just {$K$--nice\/}) whenever 
\begin{enumerate}
\item[$(\boxtimes_{\rm nice}^K)$] if so far Generic used $\st$ and $s_i$ is
  given to her as a move at a stage $i<\omega$, then  
\begin{itemize}
\item $s_i\subseteq \bigcup\limits_{j\leq i+1} {}^j (n+1)$,
  $\max(s_i)\subseteq {}^{(i+1)}(n+1)$ and
\item if $\eta\in\max(s_i)$ and $i\notin K$, then $\eta(i)=n$, and 
\item if $\eta\in\max(s_i)$ and $i\in K$, then $\suc_{s_i}(\eta\rest i)=n$,
\item if $i\in K$, $\langle \eta^i_\ell:\ell<k\rangle$ is an enumeration of
  $\max(s_i)$ and $\langle p^i_{\eta^i_\ell},q^i_{\eta^i_\ell}:\ell<k
  \rangle$ is the result of the subgame of level $i$ in which Generic uses
  $\st$, then the conditions $p^i_{\eta^i_\ell}$ (for $\ell<k$) are pairwise
  incompatible. 
\end{itemize}
\end{enumerate}
In a similar way we define when a strategy $\st$ for Generic in
$\Gsg(p,\bbP)$ or $\tGs(p,\bbP)$ is {\em nice for $K$}. 
\item We say that $\bbP$ has the {\em nice $\odot_n$--property\/} if for
every $K\in [\omega]^{\textstyle\omega}$ and $p\in\bbP$, Generic has a
$K$--nice winning strategy in $\wGs(p,\bbP)$.
\end{enumerate}
\end{definition}

\begin{remark}
\label{remgames}
\begin{enumerate}
\item At a stage $i<\lambda$ of a play of $\Gsg(p,\bbP)$, Antigeneric may play
stronger conditions, and we may require that if $\bar{q}^i=\langle
q^i_\eta:\eta\in \max(s_i)\rangle$ is his move, then the conditions 
$q^i_\eta$ are pairwise incompatible. Thus the winning criterion
$(\circledast)$ could be replaced by
\begin{enumerate}
\item[$(\circledast)^*$]\qquad there are a condition $q\geq p$ and a
$\bbP$--name $\name{\rho}$ such that 
\[q\forces_{\bbP}\mbox{`` }\name{\rho}\in \big[\bigcup_{i<\omega}s_i\big]
\ \&\ \big(\forall i<\omega\big)\big(\exists \ell<\omega\big)\big(
\name{\rho}\rest\ell\in \max(s_i)\ \&\ q^i_{\name{\rho}\rest\ell}\in
\Gamma_{\bbP}\big)\mbox{ ''.}\]
\end{enumerate}
This would make the game $\Gsg$ more like the game of \cite[Definition
A.2.1]{RoSh:777}.
\item If Generic has a winning strategy in $\Gsg(p,\bbP)$ and $K\subseteq
\omega$ is infinite, then Generic has a $K$--nice winning strategy in
$\Gsg(p,\bbP)$.   
\end{enumerate}
\end{remark}

\begin{observation}
\label{easyobs}
For a forcing notion $\bbP$ the following implications hold:
\[\begin{array}{rrl}
\mbox{$\oplus_n$--property}\ \ \Rightarrow&
\mbox{nice $\odot_n$--property}\ \Rightarrow\
\mbox{$\odot_n$--property}\ \Rightarrow &
\mbox{$\ominus_n$--property}\\
\Downarrow\qquad\qquad&&\qquad\Downarrow\\
\mbox{$n$--localization property} & & \mbox{proper.}
  \end{array}\]
\end{observation}

Let us recall definitions of forcing notions that are main examples for the
properties introduced in \ref{da2}. 

\begin{definition}
\label{forcings}
\begin{enumerate}
\item {\bf The $n$--Sacks forcing notion $\bbD_n$} consists of perfect trees 
  $p\subseteq {}^{\omega{>}}n$ such that 
\[(\forall \eta\in p)(\exists\nu\in p)(\eta\vtl\nu\ \&\ \suc_p(\nu)=n).\]
The order of $\bbD_n$ is the reverse inclusion, i.e., $p\leq_{\bbD_n} q$ if
and only if $q\subseteq p$.

\item {\bf The uniform $n$--Sacks forcing notion $\bbQ_n$} consists of
  perfect trees $p\subseteq {}^{\omega{>}}n$ such that 
\[(\exists X\in [\omega]^\omega)(\forall \eta\in p)(\lh(\eta)\in X\
  \Rightarrow\ \suc_p(\nu)=n).\]
The order of $\bbQ_n$ is the reverse inclusion, i.e., $p\leq_{\bbQ_n} q$ if
and only if $q\subseteq p$.

\item {\bf The $n$--Silver forcing notion $\bbS_n$} consists of partial
  functions $p$ such that $\Dom(p)\subseteq \omega$, $\Rng(p)\subseteq n$
  and $\omega\setminus\Dom(p)$ is infinite. The order of $\bbS_n$ is the
  inclusion, i.e., $p\leq_{\bbQ_n} q$ if and only  if $p\subseteq q$.

\item Let us assume that $G=(V,E)$ is a hypergraph on a Polish space $V$
  which is
\begin{itemize}
\item {\em $(n+1)$--regular open,\/} that is $E\subseteq [V]^{n+1}$ is open
  in the topology inherited from $V^{n+1}$, and 
\item {\em transitive,\/} that is  $\big(\forall e\in E\big)\big(\forall
v\in V\setminus e\big)\big(\exists w\in e\big)\big((e\setminus\{w\})\cup
\{v\}\in E\big)$,    
\item {\em uncountably chromatic on every open set,\/} that is for every
  non-empty open subset $U$ of $V$ and every countable family $\cF$ of
  subsets of $U$, either $\bigcup\cF\neq U$ or $[F]^{n+1}\cap
  E\neq\emptyset$ for some $F\in\cF$.  
\end{itemize}
{\bf The Geschke forcing notion $\bbP_G$ for $G$} consists of all closed
sets $C\subseteq V$ such that the hypergraph $(C,E\cap [C]^{n+1})$ is
uncountably chromatic on every non-empty open subset of $C$. The order of
$\bbP_G$ is the inverse inclusion, i.e., $C\leq_{\bbP_G} D$ if and only  if
$D\subseteq C$. 
\end{enumerate}
\end{definition}

\begin{observation}
\begin{enumerate}
\item The $n$--Sacks forcing notion $\bbD_n$ has the $\oplus_n$--property.  
\item The uniform $n$--Sacks forcing notion $\bbQ_n$ and the $n$--Silver
forcing notion $\bbS_n$ have the nice $\odot_n$--property.  
\item Assume that $G=(V,E)$ is a transitive $(n+1)$--regular open hypergraph
on a Polish space $V$ which is uncountably chromatic on every open set. Then
the corresponding Geschke forcing notion $\bbP_G$ has the
$\oplus_n$--property.  
\end{enumerate}
\end{observation}

\begin{proof}
(1)--(3)\quad Straightforward.

\noindent (4)\quad This is included in the proof of \cite[Lemma 2.8]{Ge05}. 
\end{proof}

The proofs of our theorems resemble arguments from the pre-proper era of
iterated forcing and their crucial ingredients are trees of conditions. Let
us first recall the relevant notions --- in the definition below we follow
the pattern that recently has been used in the context of iterations with
uncountable supports.   

\begin{definition}
[cf {\cite[Def. A.1.7]{RoSh:777}}, {\cite[A.3.3, A.3.2]{Sh:587}}] 
\label{dA.5}
Let $\gamma$ be an ordinal and let $\bar{\bbQ}=\langle\bbP_\xi,
\name{\bbQ}_\xi:\xi<\gamma\rangle$ be a CS iteration. 
\begin{enumerate}
\item Let $m<\omega$ and $w\subseteq \gamma$ be finite. {\em A standard
  $(w,m)^\gamma$--tree\/} is a pair $\cT=(T,\rk)$ such that  
\begin{itemize}
\item $(T,\vtl)$ is a tree with root $\langle\rangle$, $\rk:T\longrightarrow 
  w\cup\{\gamma\}$, and 
\item if $t\in T$ and $\rk(t)=\vare$, then $t$ is a sequence $\langle
  (t)_\zeta: \zeta\in w\cap\vare\rangle$, where each $(t)_\zeta$ is a
  sequence of length $m$.  
\end{itemize}
We will keep the convention that $\cT^x_y$ is $(T^x_y,\rk^x_y)$.
\item Suppose that $w_0\subseteq w_1$ are finite subsets of $\gamma$, $m_0
\leq m_1$, and $\cT_1=(T_1,\rk_1)$ is a standard $(w_1,m_1)^\gamma$--tree. 
{\em The projection $\proj^{(w_1,m_1)}_{(w_0,m_0)}(\cT_1)$ of $\cT_1$ onto
$(w_0,m_0)$} is defined as a standard $(w_0,m_0)^\gamma$--tree $\cT_0=
(T_0,\rk_0)$ such that    
\[T_0=\{\langle (t)_\zeta\rest m_0:\zeta\in w_0\cap\rk_1(t)\rangle:
t=\langle (t)_\zeta:\zeta\in w_1\cap\rk_1(t)\rangle\in T_1\}.\]
The mapping 
\[T_1\ni\langle (t)_\zeta:\zeta\in w_1\cap\rk_1(t)\rangle\longmapsto\langle  
(t)_\zeta\rest m_0: \zeta\in w_0\cap\rk_1(t)\rangle\in T_0\]
will be denoted $\proj^{(w_1,m_1)}_{(w_0,m_0)}$ too.
\item {\em A standard tree of conditions in $\bar{\bbQ}$} is a system
$\bar{p}=\langle p_t:t\in T\rangle$ such that   
\begin{itemize}
\item $(T,\rk)$ is a standard $(w,m)^\gamma$--tree for some finite set
  $w\subseteq\gamma$ and an integer $m<\omega$, 
\item $p_t\in\bbP_{\rk(t)}$ for $t\in T$, and
\item if $s,t\in T$, $s\vtl t$, then $p_s=p_t\rest\rk(s)$. 
\end{itemize}
\item Let $\bar{p}^0,\bar{p}^1$ be standard trees of conditions in
$\bar{\bbQ}$, $\bar{p}^i=\langle p^i_t:t\in T_i\rangle$, where $\cT_0=
\proj^{(w_1,m_1)}_{(w_0,m_0)}(\cT_1)$, $w_0\subseteq w_1\subseteq
\gamma$, $m_0\leq m_1$. We will write $\bar{p}^0\leq^{w_1,m_1}_{w_0,m_0}
\bar{p}^1$ (or just $\bar{p}^0\leq \bar{p}^1$) whenever for each $t\in T_1$,
letting $t'=\proj^{(w_1,m_1)}_{(w_0,m_0)}(t)\in T_0$, we have $p^0_{t'}\rest 
\rk_1(t)\leq p^1_t$. 
\end{enumerate}
\end{definition}

\begin{lemma}
\label{treelem}
Assume that 
\begin{itemize}
\item $\bar{\bbQ}=\langle\bbP_\xi,\name{\bbQ}_\xi:\xi<\gamma\rangle$ is a CS
iteration,
\item $(T,\rk)$ is a standard $(w,m)^\gamma$--tree, $w\in [\gamma]^{
\textstyle{<}\omega}$ and $\bar{p}=\langle p_t:t\in T\rangle$ is a standard
  tree of conditions in $\bar{\bbQ}$, and 
\item $\name{\tau}$ is a $\bbP_\gamma$--name for an element of $\baire$ such
  that $\forces_{\bbP_\gamma}\big(\forall\alpha<\gamma\big)\big(\name{\tau}
  \notin \bV^{\bbP_\alpha}\big)$.
\end{itemize}
Then there are a tree of conditions $\bar{q}=\langle q_t:t\in T\rangle$ and
$N\in\omega$ such that 
\begin{itemize}
\item $\bar{p}\leq\bar{q}$,
\item if $t\in T$, $\rk(t)=\gamma$, then the condition $q_t$ decides
$\name{\tau}\rest N$, say $q_t\forces_{\bbP_\gamma}\name{\tau}\rest N=
\sigma_t$, 
\item if $t_0,t_1\in T$, $\rk(t_0)=\rk(t_1)=\gamma$ and $t_0\neq t_1$, then
$\sigma_{t_0}\neq \sigma_{t_1}$. 
\end{itemize}
\end{lemma}

\begin{proof}
For $\alpha<\beta\leq\gamma$, $\name{\bbP}_{\alpha\beta}$ is a
$\bbP_\alpha$--name for a forcing notion with universe $P_{\alpha\beta}=
\{p\rest [\alpha,\beta): p\in\bbP_\beta\}$ such that 
\begin{quotation}
if $G_\alpha\subseteq\bbP_\alpha$ is generic over $\bV$ and $f,g\in
P_{\alpha\beta}$, \\
then $\bV[G_\alpha]\models f\leq_{\name{\bbP}_{\alpha\beta}[G_\alpha]} g$ if
and only if $(\exists p\in G_\alpha)(p\cup f\leq _{\bbP_\beta} p\cup g)$.  
\end{quotation}
Note that $P_{\alpha\beta}$ is from $\bV$, it is only the relation
$\leq_{\name{\bbP}_{\alpha\beta}}$ which is defined in
$\bV[G_\alpha]$. Also, $\bbP_\beta$ is isomorphic with a dense subset of
$\bbP_\alpha*\name{\bbP}_{\alpha\beta}$. 

Let $w=\{\xi_0,\ldots,\xi_i\}$ be the decreasing enumeration. We may assume
that $0\in w$ and thus $\xi_0=\max(w)$ and $\xi_i=0$. Let $M=|T|+7$,
$T^*=\{t\in T:\rk(t)=\gamma\}$. By induction on $j\leq i$ we will define 
$\bbP_{\xi_j}$--names $\name{n}^j_t$, $\name{\sigma}^j_{t,k}$ and
$\name{q}^j_{t,k}$ for $t\in T^*$and $k<M$ so that 
\begin{enumerate}
\item[(a)$_j$] $\forces_{\bbP_{\xi_j}}$`` $\name{n}^j_t\in\omega\ \&\
\name{\sigma}^j_{t,k}:\name{n}^j_t\longrightarrow\omega$ '',
\item[(b)$_j$] $\forces_{\bbP_{\xi_j}}$`` $\name{q}^j_{t,k}\in
\name{\bbP}_{\xi_j\gamma}\ \&\ p_t\rest[\xi_j,\gamma)
\leq_{\name{\bbP}_{\xi_j\gamma}}\name{q}^j_{t,k}\ \&\
\name{q}^j_{t,k}\forces_{\name{\bbP}_{\xi_j\gamma}}\name{\tau}\rest
\name{n}^j_t=\name{\sigma}^j_{t,k}$ '', 
\item[(c)$_j$] if $t\in T^*$ and $k<\ell<M$, then 
\[\forces_{\bbP_{\xi_j}}\mbox{`` }\name{q}^j_{t,k}\rest [\xi_j,\xi_0)=
\name{q}^j_{t,\ell}\rest [\xi_j,\xi_0)\ \&\ \name{\sigma}^j_{t,k} \neq
\name{\sigma}^j_{t,\ell}\mbox{ '',}\]
\item[(d)$_j$] if $t_0,t_1\in T^*$, $t=t_0\cap t_1$ and $\rk(t)>\xi_j$, then
\[\forces_{\bbP_{\xi_j}}\mbox{`` }\name{q}^j_{t_0,k}\rest \rk(t)=
\name{q}^j_{t_1,k}\rest \rk(t)\mbox{ ''.}\]
\end{enumerate}
To start the inductive process suppose that $G_{\xi_0}\subseteq\bbP_{\xi_0}$
is generic over $\bV$ and work in $\bV[G_{\xi_0}]$ for a moment. Note that
$\name{\tau}$ may be thought of as a $\name{\bbP}_{\xi_0\gamma}[G_{\xi_0}
]$--name for an element of $\baire$ such that $\forces_{
\name{\bbP}_{\xi_0\gamma}[G_{\xi_0}]}\name{\tau}\notin\bV[G_{\xi_0}]$.
Therefore we may find $n_t\in \omega$, $\sigma_{t,k}:n_t\longrightarrow
\omega$ and  $q_{t,k}\in\name{\bbP}_{\xi_0\gamma}[G_{\xi_0}]$ (for $t\in
T^*$ and $k<M$) such that for each $t\in T^*$ and $\ell,k<M$, $\ell\neq k$: 
\begin{itemize}
\item $p_t\rest [\xi_0,\gamma)\leq_{\name{\bbP}_{\xi_0\gamma}[G_{\xi_0}]}
q_{t,k}$, 
\item $q_{t,k}\forces_{\name{\bbP}_{\xi_0\gamma}[G_{\xi_0}]}\name{\tau}\rest
n_t=\sigma_{t,k}$,
\item $\sigma_{t,k}\neq\sigma_{t,\ell}$.
\end{itemize}
Now, let $\name{n}^0_t,\name{\sigma}^0_{t,k},\name{q}^0_{t,k}$ (for $t\in
T^*$, $k<M$) be $\bbP_{\xi_0}$--names for objects with properties as those
of $n_t,\sigma_{t,k},q_{t,k}$ above. 

Suppose that $j<i$ and we have defined $\bbP_{\xi_j}$--names $\name{n}^j_t,
\name{\sigma}^j_{t,k},\name{q}^j_{t,k}$ so that (a)$_j$--(d)$_j$ are
satisfied. Let $G_{\xi_{j+1}}\subseteq\bbP_{\xi_{j+1}}$ be generic over
$\bV$ and work in $\bV[G_{\xi_{j+1}}]$ for a moment. For each $s\in T$ of
rank $\rk(s)=\xi_j$ we may pick a condition $q_s\in\name{\bbP}_{
\xi_{j+1}\xi_j}[G_{\xi_{j+1}}]$ stronger than $p_s\rest [\xi_{j+1},\xi_j)$
and also we may choose $n_t\in\omega$, $\sigma_{t,k}:n_t\longrightarrow
\omega$ and $q_{t,k}$ (for $k<M$, $t\in T^*$) such that 
\[q_s\forces_{\name{\bbP}_{\xi_{j+1}\xi_j}[G_{\xi_{j+1}}]}\!\big(\forall
k\!<\! M\big)\big(\forall t\!\in\! T^*\big)\big(s\vtl t\ \Rightarrow\
[\name{n}^j_t=n_t \ \&\ \name{q}^j_{t,k}=q_{t,k}\ \&\
  \name{\sigma}^j_{t,k}=\sigma_{t,k}] \big).\] 
Now let $\name{n}^{j+1}_t, \name{\sigma}^{j+1}_{t,k}$ be
$\bbP_{\xi_{j+1}}$--names for $n_t,\sigma_{t,k}$ as above, and let
$\name{q}^{j+1}_{t,k}$ be a $\bbP_{\xi_j+1}$--name for $q_{t\rest\xi_j}\conc 
q_{t,k}$. One easily verifies that demands (a)$_{j+1}$--(d)$_{j+1}$ are
satisfied. 

Finally note that (as $\xi_i=0$) $\name{n}^i_t, \name{\sigma}^i_{t,k}$ and 
$\name{q}^i_{t,k}$ are actually objects in $\bV$, not names.

Let $\cT^+=(T^+,\rk^+)$ be a standard $(w,m+1)^\gamma$--tree such that
$\proj^{w,m+1}_{w,m}(\cT^+)=\cT$ and 
\begin{quotation}
if $t=\langle (t)_\xi:\xi\in w\rangle\in T^+$, $\rk^+(t)=\gamma$,\\
then $(t)_{\xi_0}(m)<M$ and $\big(\forall\xi\in w\cap\xi_0\big)\big(
(t)_\xi(m)=*\big)$. 
\end{quotation}
It should be clear that $\langle q^i_{t,k}:k<M\ \&\ t\in T^*\rangle$
determines a tree of conditions $\bar{q}'=\langle q_t':t'\in T^+\rangle$
such that $\bar{p}\leq \bar{q}'$ and $q_{t'}'=q^i_{t,k}$ whenever $t'\in
T^+$, $\rk^+(t')=\gamma$, $t=\proj^{w,m+1}_{w,m}(t')$ and $k=(t')_{\xi_0}
(m)$. Let $N=\max(\{n^i_t:t\in T^*\})$. Carrying out a procedure similar to
that described above we may find a tree of conditions $\bar{q}^*=\langle
q^*_s:s\in T^+\rangle$ such that $\bar{q}^*\geq \bar{q}'$ and for some
$\rho_s\in {}^n\omega$ (for $s\in T^+$, $\rk^+(s)=\gamma$) we have
\begin{itemize}
\item $q^*_s\forces_{\bbP_\gamma}\name{\tau}\rest N=\rho_s$, and 
\item if $s_0,s_1\in T^+$, $\rk^+(s_0)=\rk^+(s_1)=\gamma$, $\proj^{w,
m+1}_{w,m}(s_0)=\proj^{w,m+1}_{w,m}(s_1)$, $(s_0)_{\xi_0}(m)\neq
(t_1)_{\xi_0}(m)$, then $\rho_{t_0}\neq \rho_{t_1}$. 
\end{itemize}
Then for each $t\in T^*$ we may choose $s_t\in T^+$ such that
$\proj^{w,m+1}_{w,m}(s_t)=t$ and $\rho_{s_{t_0}}\neq \rho_{s_{t_1}}$ for
distinct $t_0,t_1\in T^*$. The choice of $\bar{q}$ should be clear now.
\end{proof}

\section{$\oplus_n$--property and CS iterations}
Here we show that CS iterations of forcing notions with $\oplus_n$--property
result in forcings with the $n$--localization property. This result covers
examples like the $n$--Sacks forcing notion $\bbD_n$ or the suitable Geschke
forcings $\bbP_G$. However, we do not know if the uniform $n$--Sacks forcing
fits here, so in the next section we will prove a result applicable to a
larger family of forcing notions. Still we believe that the proof of
\ref{iter} below is a nice preparation for the arguments in the following
section. 

\begin{theorem}
\label{iter}
Let $\bar{\bbQ}=\langle\bbP_\xi,\name{\bbQ}_\xi:\xi<\gamma\rangle$ be a
CS iteration such that for every $\xi<\gamma$,  
\[\forces_{\bbP_\xi}\mbox{`` $\name{\bbQ}_\xi$ has the $\oplus_n$--property
''.}\]  
Then 
\begin{enumerate}
\item $\bbP_{\gamma}=\lim(\bar{\bbQ})$ has the $\odot_n$--property. 
\item $\bbP_{\gamma}=\lim(\bar{\bbQ})$ has the $n$-localization
  property.    
\end{enumerate}
\end{theorem}

\begin{proof}
(1)\qquad Let $p\in\bbP_\gamma$. We are going to describe a strategy $\st$
for Generic in the game $\wGs(p,\bbP_\gamma)$. This strategy will give
Generic, at a stage $i<\omega$, a standard  $(w_i,(i+1))^\gamma$--tree
$\cT_i=(T_i,\rk_i)$. These standard trees will satisfy $\cT_i=
\proj^{(w_{i+1},i+2)}_{(w_i,i+1)}(\cT_{i+1})$ and $\{t\in T_i:\rk_i(t)=
\gamma\}$ will correspond to $\max(s_i)$ in the rules of the game. If only
we make sure that  
\begin{enumerate}
\item[$(\oplus)_0$] for each $t'\in T_i$ with $\rk_i(t)=\gamma$ we have 
\[0<\big|\big\{ t\in T_{i+1}:\proj^{(w_{i+1},i+2)}_{(w_i,i+1)}(t)=t'\big\} 
\big|\leq n,\]
\end{enumerate}
then Generic may easily build trees $s_i$ and mappings $\pi_i:\{t\in
T_i:\rk_i(t)=\gamma\}\longrightarrow s_i$ such that 
\begin{enumerate}
\item[$(\oplus)^{\rm (a)}_0$] $\Rng(\pi_i)=\max(s_i)\subseteq{}^{(i+1)}n$,
\item[$(\oplus)^{\rm (b)}_0$] $\big(\forall t_0\in T_i\big)\big(\forall t_1\in 
  T_{i+1}\big)\big(\pi_i(t_0)\vtl \pi_{i+1}(t_1)\ \Leftrightarrow\ t_0=
  \proj^{(w_{i+1},i+2)}_{(w_i,i+1)}(t_1)\big)$,  
\item[$(\oplus)^{\rm (c)}_0$] the demands of \ref{da2}(1$(\alpha)$) hold.
\end{enumerate}
Later we will even not mention the trees $s_i$ but we will work directly
with $T_i$.

As we said, in the course of the play the strategy $\st$ will instruct
Generic to choose finite sets $w_i\subseteq\gamma$ and standard
$(w_i,i+1)^\gamma$--trees $\cT_i$. She will also pick sets $K_\xi\in
[\omega]^{\textstyle\omega}$, conditions $r_i\in\bbP_\gamma$ and
$\bar{t}^i,\bar{p}^i_*,\bar{q}^i_*,k_i,i^*_\xi,\name{\st}_\xi,s_{i,\xi}, 
\name{\bar{p}}_{i,\xi},\name{\bar{q}}_{i,\xi}$. All these objects will be
constructed so that, assuming $\langle (\cT_i,\bar{t}^i,\bar{p}^i,
\bar{q}^i):i<\omega\rangle$ is the result of a play of $\wGs(p,\bbP_\gamma)$
in which Generic used $\st$ and she determined the corresponding side
objects, the following conditions are satisfied. 

\begin{enumerate}
\item[$(\oplus)_1$] $r_0(0)=p(0)$, $w_i\in [\gamma]^{\textstyle i+1}$, $w_0=
\{0\}$, $w_i\subseteq w_{i+1}$ and $\bigcup\limits_{i<\omega}\Dom(r_i)=
\bigcup\limits_{i<\omega} w_i$.   
\item[$(\oplus)_2$] If $j<i<\omega$, then $\big(\forall\xi\in w_{j+1}\big)
\big(r_j(\xi)= r_i(\xi)\big)$ and $p\leq r_j\leq r_i$. 
\item[$(\oplus)_3$] If $\xi\in w_i$, then $K_\xi\in [\omega]^{\textstyle
\omega}$ is known at stage $i$ of the play and if $\xi,\zeta\in
\bigcup\limits_{i\in\omega} w_i$ are distinct, then $K_\zeta\cap
K_\xi=\emptyset$.  
\item[$(\oplus)_4$] For $\xi\in\bigcup\limits_{i<\omega} w_i$ we have
$i^*_\xi=\min(\{i:\xi\in w_i\})\leq\min(K_\xi)$, and $\name{\st}_\xi$ is a
$\bbP_\xi$--name for a winning strategy of Generic in $\Gsg(r_{i^*_\xi}
(\xi),\name{\bbQ}_\xi)$ which is nice for $\{k\in\omega:k+i^*_\xi\in K_\xi
\}$ (see \ref{nice}, \ref{remgames}(2)). (So $\st_0$ is a $K_0$--nice
winning strategy of Generic in $\Gsg(r_0(0),\bbQ_0)$.)       
\item[$(\oplus)_5$] $\cT_i=(T_i,\rk_i)$ is a standard $(w_i,
i+1)^\gamma$--tree, $\cT_i=\proj^{(w_{i+1},i+2)}_{(w_i,i+1)}(\cT_{i+1})$.  
\item[$(\oplus)_6$] $\bar{p}^i_*=\langle p^i_{*,t}:t\in T_i\rangle$ and
$\bar{q}^i_*=\langle q^i_{*,t}:t\in T_i\rangle$ are standard trees of
conditions, $\bar{p}^i_*\leq\bar{q}^i_*\leq^{w_{i+1},i+2}_{w_i,i+1}
\bar{p}^{i+1}_*$.    
\item[$(\oplus)_7^0$] $\Dom(p^0_{*,t})=\big(\{0\}\cup\Dom(p)\big)\cap
\rk_0(t)$ for each $t\in T_0$ and $p^0_{*,t}(\xi)=p(\xi)$ for $\xi\in
\Dom(p^0_{*,t})\setminus\{0\}$, $t\in T_0$. 
\item[$(\oplus)_7^{i+1}$] For $t\in T_{i+1}$ we have $\Dom(p^{i+1}_{*,t})
=\big(\Dom(r_i)\cup w_{i+1}\big)\cap\rk_{i+1}(t)$ and $p^{i+1}_{*,t}(\xi)= 
r_i(\xi)$ for $\xi\in \Dom(p^{i+1}_{*,t})\setminus w_{i+1}$. 
\item[$(\oplus)_8$] $k_i=|\{t\in T_i:\rk_i(t)=\gamma\}|$, $\bar{t}^i=\langle
t^i_\ell:\ell<k_i\rangle$ is an enumeration of $\{t\in T_i:\rk_i(t)=
\gamma\}$, and for each $t\in T_i$ with $\rk_i(t)=\gamma$ we have
\[p^i_{*,t}\leq p^i_{t}\leq q^i_{t}\leq q^i_{*,t}.\]  
\item[$(\oplus)_9$] If $\xi\in w_i$, then $s_{i,\xi}\subseteq
\bigcup\limits_{j\leq i+1-i^*_\xi} {}^j(n+1)$ is an $n$--tree and
$\name{\bar{p}}_{i,\xi}=\langle\name{p}^\eta_{i,\xi}:\eta\in\max(s_{i,\xi})
\rangle$, $\name{\bar{q}}_{i,\xi}=\langle\name{q}^\eta_{i,\xi}:\eta\in
\max(s_{i,\xi})\rangle$ are $\bbP_\xi$--names for systems of conditions in 
$\name{\bbQ}_\xi$ (indexed by $\max(s_{i,\xi})$).    
\item[$(\oplus)_{10}$] For each  $\xi\in \bigcup\limits_{i<\omega} w_i$,  
\[\begin{array}{r}
\forces_{\bbP_\xi}\mbox{`` }\langle s_{i,\xi},\name{\bar{p}}_{i,\xi},
\name{\bar{q}}_{i,\xi}:i^*_\xi\leq i<\omega\rangle\mbox{ is a legal play of
}\Gsg(r_{i^*_\xi}(\xi),\name{\bbQ}_\xi)\\
\mbox{ in which Generic uses $\name{\st}_\xi$ ''.}
  \end{array}\]
\item[$(\oplus)_{11}$] If $t\in T_i$, $\rk_i(t)=\xi<\gamma$ (so $\xi\in
w_i$ and $i\geq i^*_\xi$), then 
\[\big\{(s)_\xi:t\vartriangleleft s\in T_i\big\}=\big\{\eta: \eta\rest
i^*_\xi\in {}^{i^*_\xi}\{*\}\ \ \&\ \ \big(\exists\nu\in\max(s_{i,\xi})
\big)\big(\eta=(\eta\rest i^*_\xi)\conc\nu\big)\big\}.\]  
\item[$(\oplus)_{12}$] If $t\in T_i$, $\xi<\rk_i(t)$, $\xi\in w_i$ and
$(t)_\xi=\big((t)_\xi\rest i^*_\xi\big)\conc\nu$ (so $\nu\in\max(s_{i,
\xi})$), then 
\[p^i_{*,t}\rest\xi\forces_{\bbP_\xi}\mbox{`` }p^i_{*,t}(\xi)=
\name{p}^\nu_{i,\xi}\mbox{ ''\qquad and\qquad }q^i_{*,t}\rest\xi
\forces_{\bbP_\xi}\mbox{`` }q^i_{*,t}(\xi)=\name{q}^\nu_{i,\xi}\mbox{ ''.}\]     
\item[$(\oplus)_{13}$] If $t_0,t_1\in T_i$, $\rk_i(t_0)=\rk_i(t_1)$ and
$\xi\in w_i\cap\rk_i(t_0)$, $t_0\rest\xi=t_1\rest\xi$ but $\big(t_0\big)_\xi
\neq\big(t_1\big)_\xi$, then 
\[q^i_{*,t_0\rest\xi}\forces_{\bbP_\xi}\mbox{`` the conditions $q^i_{*,t_0}
(\xi),q^i_{*,t_1}(\xi)$ are incompatible ''.}\]  
\item[$(\oplus)_{14}$] $\Dom(r_i)=\bigcup\limits_{t\in T_i}\Dom(q^i_{*,t})
\cup\Dom(p)$ and if $t\in T_i$, $\xi\in\Dom(r_i)\cap\rk_i(t)\setminus w_i$, 
then $q^i_{*,t}\rest\xi\forces_{\bbP_\xi}\mbox{`` }r_i(\xi)\geq
q^i_{*,t}(\xi)\mbox{ ''}$.
\end{enumerate}

To describe the instructions given by $\st$ at stage $i<\omega$ of a play of
$\wGs(p,\bbP_\gamma)$ let us assume that $\big\langle(\cT_j,\bar{t}^j,
\bar{p}^j,\bar{q}^j):j<i\big\rangle$ is the result of the play so far and
that Generic constructed aside the objects appearing in
$(\oplus)_1$--$(\oplus)_{14}$ (and they have the respective properties). 

For definiteness of our definitions, whenever we say ``Generic chooses/picks
$X$ such that'' we really mean ``Generic takes the $<^*_\chi$--first $X$
such that''.

First, Generic uses her favourite bookkeeping device to determine $w_i$ such
that the demands in $(\oplus)_1$ are satisfied (and that at the end we will
have $\bigcup\limits_{j<\omega}\Dom(r_j)=\bigcup\limits_{j<\omega}w_j$) and
then again she uses the bookkeeping device to determine $K_\xi$ so that
$(\oplus)_3+(\oplus)_4$ hold. Note that $i^*_\xi$ for $\xi\in w_i$ is
defined by $(\oplus)_4$, also the choice of $\name{\st}_\xi$ is determined
by $(\oplus)_4$ (remember that $r_{i^*_\xi}(\xi)$ is determined by
$(\oplus)_2$). 

Now $(\oplus)_9+(\oplus)_{10}$ decide $s_{i,\xi}$ (for $\xi\in w_i$) 
and since $\name{\st}_\xi$ is (a name for) a nice for $K_\xi-i^*_\xi$
strategy, we know that $s_{i,\xi}$ can be easily read from the truth value
of ``$i+i^*_\xi\in K_\xi$''. Plainly $\max(s_{i,\xi})\subseteq {}^{i+1-
i^*_\xi}(n+1)$ and the clauses mentioned before determine
$\name{\bar{p}}_{i,\xi}=\langle\name{p}^\eta_{i,\xi}:\eta\in\max(s_{i,\xi})
\rangle$.  Now the choice of the standard tree $\cT_i$ is fully described by
$(\oplus)_5+(\oplus)_{11}$ and clearly $(\oplus)_0$ holds then too. For each
$t\in T_i$ Generic picks a condition $p^i_{*,t}\in\bbP_{\rk_i(t)}$ so that
the demands of $(\oplus)_7^i+(\oplus)_{12}$ are satisfied. (One may use
$(\oplus)_{14}$ to argue that the last demand in $(\oplus)_6$ is satisfied.)

After the above choices are made, Generic (in the play of $\wGs(p,
\bbP_\gamma)$) puts $T_i$ as her inning and  Antigeneric chooses an
enumeration $\bar{t}^i=\langle t^i_\ell:\ell<k_i\rangle$ of $\{t\in
T_i:\rk_i(t)=\gamma\}$. Now the two players start a subgame of length
$k_i=|\{t\in T_i:\rk_i(t)= \gamma\}|$. During the subgame Generic will also
pick (for temporary use) trees of conditions $\bar{q}^\temp_\ell=\langle
q^\temp_{\ell,t}:t\in T_i \rangle$ for $\ell\leq k_i$. So, she lets
$\bar{q}^\temp_0=\bar{p}^i_*$ and she plays $p^i_{t^i_0}=q^\temp_{0,t^i_0}$
as her first inning in the subgame.  Antigeneric answers with
$q^i_{t^i_0}\geq p^i_{t^i_0}$ after which Generic picks a tree of conditions
$\bar{q}^\temp_1$ so that $q^\temp_{1,t^i_0}=q^i_{t^i_0}$ and for each
$0<\ell<k_i$ 
\begin{itemize}
\item if $t\vtl t^i_0$, $t\vtl t^i_\ell$ and $\rk_i(t)$ is the largest
  possible, then 
\[q^\temp_{1,t^i_\ell}=q^i_{t^i_0}\rest\rk_i(t)\conc q^\temp_{0,t^i_\ell}
  \rest [\rk_i(t),\gamma).\]
\end{itemize}
Now, if the players arrived to level $\ell^*<k_i$ of the subgame and
$\bar{q}^\temp_{\ell^*}$ was chosen, then Generic plays
$p^i_{t^i_{\ell^*}}=q^\temp_{\ell^*,t^i_{\ell^*}}$. After Antigeneric
answered with $q^i_{t^i_{\ell^*}}\geq p^i_{t^i_{\ell^*}}$, Generic builds a
tree of conditions $\bar{q}^\temp_{\ell^*+1}$ so that $q^\temp_{\ell^*+1,
t^i_{\ell^*}}=q^i_{t^i_{\ell^*}}$ and for each $\ell<k_i$, $\ell\neq\ell^*$
\begin{itemize} 
\item if $t\vtl t^i_{\ell^*}$, $t\vtl t^i_\ell$ and $\rk_i(t)$ is the
  largest possible, then 
\[q^\temp_{\ell^*+1,t^i_\ell}=q^i_{t^i_{\ell^*}}\rest\rk_i(t)\conc
  q^\temp_{\ell^*,t^i_\ell}\rest [\rk_i(t),\gamma).\]
\end{itemize}

When the subgame is over Generic lets $\bar{q}^i_*=\bar{q}^\temp_{k_i}$. Note
that the demand of $(\oplus)_{13}$ is satisfied because the strategies 
$\st_\xi$ are nice, also the relevant parts of $(\oplus)_6+(\oplus)_8$
hold. The names $\name{\bar{q}}_{i,\xi}$ (for $\xi\in w_i$) are chosen so
that $q^i_{*,t}\rest\xi\forces_{\bbP_\xi}$``$\name{q}^\nu_{i,\xi}=
q^i_{*,t}(\xi)$'' and $\forces_{\bbP_\xi}$``$\name{p}^\nu_{i,\xi}\leq
\name{q}^\nu_{i,\xi}$'' whenever $t\in T_i$, $\xi\in w_i\cap\rk_i(t)$,
$\big((t)_\xi\rest i^*_\xi \big)\conc\nu=(t)_\xi$, $\nu\in\max(s_{i,\xi})$.
Then $(\oplus)_9+(\oplus)_{12}$ are satisfied. Finally Generic chooses
$r_i\in\bbP_\gamma$ essentially by conditions $(\oplus)_2+(\oplus)_{14}$
(and our rule of picking ``the $<^*_\chi$--first such that'').   

This completes the description of the side objects constructed by Generic
and her innings at stage $i$. We also verified that clauses
$(\oplus)_0$--$(\oplus)_{14}$ hold and thus the description of the strategy
is complete.   
\medskip

We are going to argue that $\st$ is a winning strategy for Generic in
$\wGs(p,\bbP_\gamma)$. To this end suppose that $\langle (\cT_i,\bar{t}^i,
\bar{p}^i,\bar{q}^i):i<\omega\rangle$ is the result of a play of
$\wGs(p,\bbP_\gamma)$ in which Generic used $\st$, and the objects
constructed at each stage $i<\omega$ are 
\begin{enumerate}
\item[$(\boxdot)$] $w_i,\cT_i,,\bar{t}^i,\bar{p}^i,\bar{q}^i,r_i,
  \bar{p}^i_*,\bar{q}^i_*,k_i,K_\xi,i^*_\xi, \name{\st}_\xi,s_{i,\xi},
  \name{\bar{p}}_{i,\xi},\name{\bar{q}}_{i,\xi}$ for $\xi\in w_i$, 
\end{enumerate}
and they satisfy the requirements $(\oplus)_0$---$(\oplus)_{14}$. 

We define a condition $q\in\bbP_\gamma$ as follows. Let $\Dom(q)=
\bigcup\limits_{i<\omega} w_i=\bigcup\limits_{i<\omega}\Dom(r_i)$ and for
$\xi\in\Dom(q)$ let $q(\xi)$ be a $\bbP_\xi$--name for a condition in
$\name{\bbQ}_\xi$ such that 
\[\forces_{\bbP_\xi}\mbox{`` }q(\xi)\geq r_{i^*_\xi}(\xi)\mbox{ and } q(\xi)
\forces_{\name{\bbQ}_\xi} \big(\forall i\geq i^*_\xi\big)\big(\exists \nu\in
\max(s_{i,\xi})\big)\big(\name{q}^\nu_{i,\xi}\in\Gamma_{\name{\bbQ}_\xi}\big)
\mbox{ ''.}\]
Clearly $q$ is well defined (remember $(\oplus)_{10}$) and $q\geq p$ (remember
$(\oplus)_1+(\oplus)_2$). Also $q\geq r_i$ for all $i<\omega$. 

We will show that for each $i<\omega$ the family $\{q^i_{*,t}: t\in T_i\ \&\
\rk_i(t)=\gamma\}$ is predense above $q$ (and this clearly will imply that
Generic won the play). So suppose $q^+\geq q$, $i<\omega$ and $w_i\cup\{
\gamma\}=\{\xi_0,\xi_1,\ldots,\xi_i,\xi_{i+1}\}$ (the increasing
enumeration, so $\xi_0=0$). By induction on $j\leq i$ we choose an
increasing sequence $\langle q^+_j:j\leq i\rangle\subseteq\bbP_\gamma$ and
we will also define $t\rest \xi_j+1$. 

First, by the choice of $q(0)$ there is $\nu\in\max(s_{i,\xi_0})$ such that
the conditions $q^+(0)$ and $q^\nu_{i,\xi_0}=\name{q}^\nu_{i,\xi_0}$ are
compatible. Let $(t)_0=\nu$ (this defines $t\rest(\xi_0+1)$). Let $q^+_0\in
\bbP_{\xi_1}$ be such that $q^+_0(0)$ is stronger than both $q^+(0)$ and
$q^\nu_{i,\xi_0}$, and let $q^+_0\rest(\xi_0,\xi_1)=q^+\rest(\xi_0,\xi_1)$.
It follows from $(\oplus)_{12}+(\oplus)_{14}$ that  $q^+_0$ is stronger than
$q^i_{*,t\rest(\xi_0+1)}$ (and, of course, it is stronger than $q^+\rest
\xi_1$). Now suppose that $j<i$ and we have defined $t\rest(\xi_j+1)\in T_i$
and a condition $q^+_j\in \bbP_{\xi_{j+1}}$ stronger than both $q^+\rest
\xi_{j+1}$ and $q^i_{*,t\rest(\xi_j+1)}$. Necessarily  
\[q^+_j\forces_{\bbP_{\xi_{j+1}}}\mbox{`` }\big(\exists\nu\in
\max(s_{i,\xi_{j+1}})\big)\big(\name{q}^\nu_{i,\xi_{j+1}}, q^+(\xi_{j+1})
\mbox{ are compatible }\big)\mbox{ ''}\]
so we may choose $\nu\in\max(s_{i,\xi_{j+1}})$ and a condition $q_{j+1}\in
\bbP_{\xi_{j+1}}$ stronger than $q^+_j$ such that  
\[q_{j+1}\forces_{\bbP_{\xi_{j+1}}}\mbox{`` }\name{q}^\nu_{i,\xi_{j+1}},
q^+(\xi_{j+1})\mbox{ are compatible ''}.\] 
Let $(t)_{\xi_{j+1}}=\langle *\ldots*\rangle\conc\nu$ (thus $t\rest(\xi_{j+
1}+1)$ has been defined) and let $q^+_{j+1}\in \bbP_{\xi_{j+2}}$ be such
that $q^+_{j+1}\rest\xi_{j+1}=q_{j+1}$,  
\[q^+_{j+1}\rest\xi_{j+1}\forces_{\bbP_{\xi_{j+1}}}\mbox{`` }q^+_{j+1}(
\xi_{j+1})\geq \name{q}^\nu_{i,\xi_{j+1}}\quad \&\quad q^+_{j+1}(\xi_{j+1})
\geq q^+(\xi_{j+1})\mbox{ ''}\]
and $q^+_{j+1}\rest(\xi_{j+1},\xi_{j+2})=q^+\rest(\xi_{j+1},\xi_{j+2})$.
Then by $(\oplus)_{12}+(\oplus)_{14}$ the condition $q^+_{j+1}$ is stronger
than $q^i_{*,t\rest(\xi_{j+1}+1)}$ and $q^+\rest\xi_{j+2}$.    

Finally look at $t=t\rest\xi_{i+1}$ and $q^+_{i+1}$. 
\medskip

\noindent (2)\qquad Since we do not know if ``the $\odot_n$--property''
implies ``the $n$--localization property'', we cannot just say that the
statement in (2) follows from (1). However, the reason for the weaker
``$\odot_n$'' in the conclusion of \ref{iter}(1) (and not ``$\oplus_n$'') is
that in our description of the strategy $\st$, we have to make sure that the
conditions played by Antigeneric form a tree of conditions. 

So to show that $\bbP_\gamma$ has the $n$--localization property we use
\ref{treelem} and the procedure described in the proof of 
\ref{iter}(1). Suppose that $\name{\tau}$ is a $\bbP_\gamma$--name for an
element of $\baire$; we may assume that $\forces_{\bbP_\gamma}\big(\forall
\alpha<\gamma\big)\big(\name{\tau}\notin\bV^{\bbP_\alpha}\big)$. Let
$p\in\bbP_\gamma$. Construct a sequence
\[\big\langle w_i,\cT_i,\bar{p}^i_*,\bar{q}^i_*,r_i,\langle i^*_\xi,
\name{\st}_\xi, K_\xi:\xi\in w_i\rangle,\bar{\sigma}^i,\bar{m}:i<\omega
\big\rangle\] 
such that conditions $(\oplus)_0$--$(\oplus)_7$ and $(\oplus)_9$--$(
\oplus)_{14}$ are satisfied and 
\begin{enumerate}
\item[$(\oplus)_{15}$] $\bar{m}=\langle m_i:i<\omega\rangle\subseteq\omega$,
$0=m_0\leq m_i<m_{i+1}$ for $i<\omega$, 
\item[$(\oplus)_{16}$] $\bar{\sigma}^i=\langle\sigma^i_t:t\in T_i\ \&\
\rk_i(t)=\gamma\rangle\subseteq {}^{[m_i,m_{i+1})}\omega$ and if $t,t'\in
T'$, $t\neq t'$, $\rk_i(t)=\rk_i(t')=\gamma$, then $\sigma^i_t\neq
\sigma^i_{t'}$, and 
\item[$(\oplus)_{17}$] $q^i_{*,t}\forces_{\bbP_\gamma}\sigma^i_t\vtl
\name{\tau}$ for $t\in T_i$, $\rk_i(t)=\gamma$.
\end{enumerate}
Then pick $q\in \bbP_\gamma$ stronger than $p$ and such that for each
$i<\omega$ the family $\{q^i_{*,t}:t\in T_i\ \&\ \rk_i(t)=\gamma\}$ is
predense above $q$ (this is done exactly as in part (1)). Let $S\subseteq
{}^{\omega{>}}\omega$ be a tree such that 
\[[S]=\big\{f\in\baire:\big(\forall i<\omega\big)\big(\exists t\in T_i
\big)\big(\rk_i(t)=\gamma\ \&\ f\rest [m_i,m_{i+1})=\sigma^i_t\big)
\big\}.\]
Then $S$ is an $n$--ary tree and $q\forces_{\bbP_\gamma}\name{\tau}\in
[S]$. 
\end{proof}

\section{$\odot_n$--property and CS iterations}
The result of the revious section is not applicable to $\bbQ_n,\bbS_n$ as
these forcing have nice $\odot_n$--property only. An iteration theorem
suitable for that property is presented below. It is not sufficient for
claiming the $n$--localization property, so later we formulate yet another
property and we argue that it implies the $n$--localization of the limits of 
CS iterations. 

\begin{theorem}
\label{seciter}
If $\bar{\bbQ}=\langle\bbP_\xi,\name{\bbQ}_\xi:\xi<\gamma\rangle$ is a
CS iteration such that for every $\xi<\gamma$,  
\[\forces_{\bbP_\xi}\mbox{`` $\name{\bbQ}_\xi$ has the nice
$\odot_n$--property '',}\]  
then $\bbP_{\gamma}=\lim(\bar{\bbQ})$ has the $\ominus_n$--property.
\end{theorem}

\begin{proof}
Let $p\in\bbP_\gamma$. We are going to describe a strategy $\st$ for Generic 
in the game $\tGs(p,\bbP_\gamma)$.  As in the proof of Theorem \ref{iter},
the strategy $\st$ will give Generic, at a stage $i<\omega$, a standard
$(w_i, (i+1))^\gamma$--tree $\cT_i=(T_i,\rk_i)$ such that $\cT_i=
\proj^{(w_{i+1},i+2)}_{(w_i,i+1)}(\cT_{i+1})$ and 
\begin{enumerate}
\item[$(\odot)_0$] for each $t'\in T_i$ with $\rk_i(t)=\gamma$ we have 
\[0<\big|\big\{ t\in T_{i+1}:\proj^{(w_{i+1},i+2)}_{(w_i,i+1)}(t)=t'\big\} 
\big|\leq n.\]
\end{enumerate}
Generic will also pick sets $K_\xi\in [\omega]^{\textstyle\omega}$,
conditions $r_i\in\bbP_\gamma$ and $k_i,i^*_\xi,\name{\st}_\xi,s_{i,\xi},
\name{\bar{p}}_{i,\xi},\name{\bar{q}}_{i,\xi}$. All these objects will be
constructed so that, assuming $\langle (\cT_i,\bar{p}^i,\bar{q}^i):i<\omega
\rangle$ is the result of a play in which Generic used $\st$ and she
determined the corresponding side objects, the following demands are
satisfied.    

\begin{enumerate}
\item[$(\odot)_1$] $r_0(0)=p(0)$, $w_i\in [\gamma]^{\textstyle i+1}$, $w_0=
\{0\}$, $w_i\subseteq w_{i+1}$ and $\bigcup\limits_{i<\omega}\Dom(r_i)=
\bigcup\limits_{i<\omega} w_i$.   
\item[$(\odot)_2$] If $j<i<\omega$, then $\big(\forall\xi\in w_{j+1}\big)
\big(r_j(\xi)= r_i(\xi)\big)$ and $p\leq r_j\leq r_i$. 
\item[$(\odot)_3$] If $\xi\in w_i$, then $K_\xi\in
[\omega]^{\textstyle\omega}$ is known at stage $i$ of the play and if
$\xi,\zeta\in\bigcup\limits_{i\in\omega} w_i$ are distinct, then $K_\zeta
\cap K_\xi=\emptyset$.   
\item[$(\odot)_4$] For $\xi\in\bigcup\limits_{i<\omega} w_i$ we have
$i^*_\xi=\min(\{i:\xi\in w_i\})\leq\min(K_\xi)$, and $\name{\st}_\xi$ is a 
$\bbP_\xi$--name for a winning strategy of Generic in $\wGs(r_{i^*_\xi}
(\xi),\name{\bbQ}_\xi)$ which is nice for $\{k\in\omega:k+i^*_\xi\in K_\xi
\}$ (see \ref{nice}). (So $\st_0$ is a $K_0$--nice winning strategy of 
Generic in $\wGs(r_0(0),\bbQ_0)$.)      
\item[$(\odot)_5$] If $\xi\in w_i$, then $s_{i,\xi}\subseteq
\bigcup\limits_{j\leq i+1-i^*_\xi} {}^j(n+1)$ is an $n$--tree and
$\name{\bar{p}}_{i,\xi}=\langle\name{p}^\eta_{i,\xi}:\eta\in\max(s_{i,\xi})
\rangle$, $\name{\bar{q}}_{i,\xi}=\langle\name{q}^\eta_{i,\xi}:\eta\in
\max(s_{i,\xi})\rangle$ are $\bbP_\xi$--names for systems of conditions in 
$\name{\bbQ}_\xi$ (indexed by $\max(s_{i,\xi})$).    
\item[$(\odot)_6$] For each  $\xi\in \bigcup\limits_{i<\omega} w_i$,  
\[\begin{array}{r}
\forces_{\bbP_\xi}\mbox{`` }\langle s_{i,\xi},\name{\bar{p}}_{i,\xi},
\name{\bar{q}}_{i,\xi}:i^*_\xi\leq i<\omega\rangle\mbox{ is a legal play of 
}\wGs(r_{i^*_\xi}(\xi),\name{\bbQ}_\xi)\\
\mbox{ in which Generic uses $\name{\st}_\xi$ and the orders of }
\max(s_{i,\xi})\\
\mbox{chosen by Antigeneric are given by $<^*_\chi$ ''.} 
\end{array}\]
\item[$(\odot)_7$] $\cT_i=(T_i,\rk_i)$ is a standard $(w_i,
i+1)^\gamma$--tree, $\cT_i=\proj^{(w_{i+1},i+2)}_{(w_i,i+1)}(\cT_{i+1})$.  
\item[$(\odot)_8$] If $t\in T_i$, $\rk_i(t)=\xi<\gamma$ (so $\xi\in w_i$ and
$i\geq i^*_\xi$), then  
\[\big\{(s)_\xi:t\vartriangleleft s\in T_i\big\}=\big\{\eta: \eta\rest
i^*_\xi\in {}^{i^*_\xi}\{*\}\ \ \&\ \ \big(\exists\nu\in\max(s_{i,\xi})
\big)\big(\eta=(\eta\rest i^*_\xi)\conc\nu\big)\big\}.\]  
\item[$(\odot)_9$] $k_i=|\{t\in T_i:\rk_i(t)=\gamma\}|$.
\item[$(\odot)_{10}$] If $\langle t^i_\ell:\ell<k_i\rangle$ is the list of
$\{t\in T_i:\rk_i(t)=\gamma\}$ chosen by Generic, $\ell<m<k_i$, $\xi\in w_i$ 
and $t_\ell^i\rest\xi=t_m^i\rest\xi$ but $\big(t_\ell^i\big)_\xi\neq
\big(t_m^i\big)_\xi$, then $q^i_{t_\ell^i}\rest\xi\leq p^i_{t_m^i}\rest
\xi$ and 
\[p^i_{t_m^i}\rest\xi\forces_{\bbP_\xi}\mbox{`` the conditions
$p^i_{t^i_\ell}(\xi),p^i_{t^i_m}(\xi)$ are incompatible ''.}\]  
\item[$(\odot)_{11}$] If $t\in T_i$, $\rk_i(t)=\gamma$, $\xi\in w_i$ and
$(t)_\xi=(t)_\xi\rest i^*_\xi\conc \eta$, $\eta\in\max(s_{i,\xi})$, then  
\[p^i_t\rest\xi\forces_{\bbP_\xi}\name{p}^\eta_{i,\xi}\leq p^i_t(\xi)\quad 
\mbox{ and }\quad q^i_t\rest\xi\forces_{\bbP_\xi}q^i_t(\xi)\leq
\name{q}^\eta_{i,\xi}.\]
\item[$(\odot)_{12}$] $\Dom(r_i)=\bigcup\limits_{t\in T_i}\Dom(q^i_t)\cup
\Dom(p)$ and if $t\in T_i$, $\rk_i(t)=\gamma$, $\xi\in\Dom(r_i)\setminus
w_i$, then $q^i_t\rest\xi\forces_{\bbP_\xi}\mbox{`` }q^i_t(\xi)\leq
r_i(\xi)$ ''.
\end{enumerate}

To describe the instructions given by $\st$ at stage $i<\omega$ of a play of
$\tGs(p,\bbP_\gamma)$ let us assume that $\big\langle(\cT_j,\bar{p}^j,
\bar{q}^j):j<i\big\rangle$ is the result of the play so far and that Generic
constructed aside the objects appearing in $(\odot)_1$--$(\odot)_{12}$ (and
they have the respective properties).  

For definiteness of our definitions, whenever we say ``Generic chooses/picks
$X$ such that'' we really mean ``Generic takes the $<^*_\chi$--first $X$
such that''.

First, Generic uses her favourite bookkeeping device to determine $w_i$ and
$K_\xi$ so that $(\odot)_1+(\odot)_3+(\odot)_4$ hold. Note that $(\odot)_4$
determines $i^*_\xi$ and $\name{\st}_\xi$ for $\xi\in w_i$ (remember that
$r_{i^*_\xi}(\xi)$ is given by $(\odot)_2$). Now $(\odot)_5+(\odot)_6$ and
the truth value of ``$i+i^*_\xi\in K_\xi$'' decide $s_{i,\xi}$ (for $\xi\in
w_i$), remember that $\name{\st}_\xi$ is (a name for) a strategy which is
nice for $K_\xi-i^*_\xi$. Plainly $\max(s_{i,\xi})\subseteq {}^{i+1-i^*_\xi}
(n+1)$. The choice of the standard tree $\cT_i$ is fully described by
$(\odot)_7+(\odot)_8$ and clearly $(\odot)_0$ holds then too. Also $k_i$ is
given by $(\odot)_9$.  

Let $\{\zeta_j:j<j^*\}$ be the increasing enumeration of $\{\xi\in w_i:
|\max(s_{i,\xi})|>1\}$ (so $j^*\leq i+1$ and we may assume that $j^*\neq
0$). We will think of $\max(s_{i,\xi})$ (for $\xi\in w_i$) as linearly 
ordered by $<^*_\chi$ (restricted suitably). This linear order determines a
list of $\max(s_{i,\xi})$ which will be considered as an inning of
Antigeneric in answer to the choice of $s_{i,\xi}$ by Generic in
$\wGs(r_{i^*_\xi}(\xi),\name{\bbQ}_\xi)$. We may identify $\{(t)_\xi:t\in
T_i\ \&\ \rk_i(t)=\xi\}$ with $\max(s_{i,\xi})$ by the mapping
$(t)_\xi\mapsto (t)_\xi\rest [i^*_\xi,i)$, so in particular the linear order
of $\max(s_{i,\xi})$ determines a linear ordering of $\{(t)_\xi:t\in T_i\
\&\ \rk_i(t)>\xi\}$ (for $\xi\in w_i$). Generic takes the lexicographic
product of these orderings for all coordinates $\xi\in w_i$ and she lets
$\langle t^i_m:m<k_i\rangle$ be the increasing (in this lexicographic order)
enumeration of $\{t\in T_i:\rk_i(t)=\gamma\}$. Then for each $t\in T_i$ with
$\rk_i(t)=\gamma$ and $\xi\in w_i$ we have
\begin{enumerate}
\item[$(\odot)_{13}$] for some $m_0<m_1\leq k_i$, 
\[\big\{t'\in T_i:\rk_i(t')=\gamma\ \&\ t'\rest\xi=t\rest\xi\big\}=\big\{
t^i_m: m_0\leq m<m_1\big\}.\]  
\end{enumerate}
The interval $[m_0,m_1)$ as in $(\odot)_{13}$ will be called {\em the 
neighbourhood of $t$ at $\xi$}.  Note that if $\xi=\zeta_j$ for some $j<j^*$,
then $m_1>m_0+1$, otherwise $m_1=m_0+1$. For $j<j^*$ let $m_j$ be such that
$[0,m_j)$ is the neighbourhood of $t^i_0$ at $\zeta_j$ (so $m_{j^*-1}<
m_{j^*-2}<\ldots<m_0=k_i$). 

Now the two players start a subgame of length $k_i$. 

For $j<m_{j^*-1}$ Generic proceeds in the subgame as follows. First
$p^i_{t^i_0}\in\bbP_\gamma$ is such that  
\[\Dom\big(p^i_{t^i_0}\big)=\left\{\begin{array}{ll}
w_0\cup\Dom(p)      &\mbox{if }i=0,\\
w_i\cup\Dom(r_{i-1})&\mbox{if }i>0,
				   \end{array}\right.\]
and if $\xi\in \Dom(p^i_{t^i_0})\setminus w_i$ then 
\[p^i_{t^i_0}(\xi)=\left\{\begin{array}{ll}
p(\xi)       &\mbox{if }i=0,\\
r_{i-1}(\xi) &\mbox{if }i>0.
			  \end{array}\right.\]
For $\xi\in w_i$, $p^i_{t^i_0}(\xi)$ is a $\bbP_\xi$--name for an element of
$\name{\bbQ}_\xi$ such that  
\[\begin{array}{ll}
\forces_{\bbP_\xi}&\mbox{`` $p^i_{t^i_0}(\xi)$ is the condition given to
Generic by $\name{\st}_\xi$ in the subgame of}\\
&\ \ \wGs(r_{i^*_\xi}(\xi),\name{\bbQ}_{\xi})\mbox{ after }\langle s_{j,\xi},
\name{\bar{p}}_{j,\xi}, \name{\bar{q}}_{j,\xi}:i^*_\xi\leq j<i\rangle \mbox{
was played and}\\
&\ \mbox{ Antigeneric picked the $<^*_\chi$--increasing enumeration of
$\max(s_{i,\xi})$ ''}  
  \end{array}\]
A straightforward induction (using $(\odot)_{12}+(\odot)_{11}$) shows that 
\begin{enumerate}
\item[$(\odot)_{14}^0$] if $i>0$, $t'=\proj^{(w_i,i+1)}_{(w_{i-1},i)}(t^i_0)$,
  then $p^i_{t^i_0}\geq q^{i-1}_{t'}$ (and $p^0_{t^0_0}\geq p$). 
\end{enumerate}
The condition $p^i_{t^i_0}$ is the first inning of Generic in the subgame,
after which Antigeneric answers with $q^i_{t^i_0}\geq p^i_{t^i_0}$. 

Now suppose that the two players arrived to a step $0<m<m_{j^*-1}$ of the
subgame. The inning of Generic now is defined similarly to that at stage 0: 
$p^i_{t^i_m}\in\bbP_\gamma$ is such that  
\[\Dom\big(p^i_{t^i_m}\big)=\left\{\begin{array}{ll}
w_0\cup\Dom(p)& \mbox{if }i=0,\\ 
w_i\cup\Dom(r_{i-1})\cup\big(\Dom(q^i_{t^i_{m-1}})\cap\zeta_{j^*-1}\big)& 
\mbox{if }i>0,
\end{array}\right.\]
and $p^i_{t^i_m}\rest\zeta_{j^*-1}=q^i_{t^i_{m-1}}\rest\zeta_{j^*-1}$, and 
$p^i_{t^i_m}(\xi)=p^i_{t^i_0}(\xi)$ for $\xi\in\Dom(p^i_{t^i_m})\setminus
(\zeta_{j^*-1}+1)$, and finally $p^i_{t^i_m}(\zeta_{j^*-1})$ is a
$\bbP_{\zeta_{j^*-1}}$--name for an element of $\name{\bbQ}_{\zeta_{j^*-1}}$
such that   
\[\begin{array}{l}
p^i_{t^i_m}\rest\zeta_{j^*-1}\forces_{\bbP_{\zeta_{j^*-1}}}\mbox{``
$p^i_{t^i_m}(\zeta_{j^*-1})$ is the inning of Generic given by }
\name{\st}_{\zeta_{j^*-1}}\mbox{ in}\\
\qquad \mbox{the subgame of level $i$ of }\wGs(r_{i^*_{\zeta_{j^*-1}}}(
\zeta_{j^*-1}),\name{\bbQ}_{\zeta_{j^*-1}})\mbox{ after the two}\\
\qquad \mbox{players played }\langle p^i_{t^i_0}(\zeta_{j^*-1}),
q^i_{t^i_0}(\zeta_{j^*-1})\rangle,\ldots,\langle p^i_{t^i_{m-1}}(
\zeta_{j^*-1}),q^i_{t^i_{m-1}}(\zeta_{j^*-1})\rangle\\
\qquad \mbox{as the conditions attached to }(t^i_0)_{\zeta_{j^*-1}},\ldots,
(t^i_{m-1})_{\zeta_{j^*-1}},\mbox{ respectively ''.}  
  \end{array}\]
Again, one may verify by induction on $\xi\in\Dom(p^i_{t^i_m})$ that 
\begin{enumerate}
\item[$(\odot)_{14}^m$] if $i>0$, $t'=\proj^{(w_i,i+1)}_{(w_{i-1},i)}(t^i_m)$,
  then $p^i_{t^i_m}\geq q^{i-1}_{t'}$ (and $p^0_{t^0_m}\geq p$). 
\end{enumerate}
Now, $p^i_{t^i_m}$ is Generic's inning at this stage of the subgame, and
after this Antigeneric answers with $q^i_{t^i_m}\geq p^i_{t^i_m}$. 

Thus we have described how Generic plays in the first $m_{j^*-1}$ steps of the
subgame ---  let us call this procedure $\proc_{j^*-1}(0,m_{j^*-1},t^i_0,
p^i_{t^i_0})$.   

Suppose that $m_{j^*-1}<k_i$ and let $m'>m_{j^*-1}$ be such that
$[m_{j^*-1},m')$ is the neighbourhood of $t^i_{m_{j^*-1}}$ at
$\zeta_{j^*-1}$ (so $j^*>1$ and so $i>0$). Let $p^i_{t^i_{m_{j^*-1}}}\in
\bbP_\gamma$ be such that   
\[\Dom\big(p^i_{t^i_{m_{j^*-1}}}\big)=w_i\cup\Dom(r_{i-1})\cup\big(
\Dom(q^i_{t^i_{m_{j^*-1}-1}})\cap\zeta_{j^*-2}\big)\]
and $p^i_{t^i_{m_{j^*-1}}}\rest\zeta_{j^*-2}=q^i_{t^i_{m_{j^*-1}-1}}\rest
\zeta_{j^*-2}$, and $p^i_{t^i_{m_{j^*-1}}}(\zeta_{j^*-2})$ is a $\bbP_{
\zeta_{j^*-2}}$--name for an element of $\name{\bbQ}_{\zeta_{j^*-2}}$ such
that   
\[\begin{array}{r}
p^i_{t^i_{m_{j^*-1}}}\rest\zeta_{j^*-2}\forces_{\bbP_{\zeta_{j^*-2}}}
\mbox{`` $p^i_{t^i_{m_{j^*-1}}}(\zeta_{j^*-2})$ is the inning of Generic
given by }\name{\st}_{\zeta_{j^*-2}}\ \\ 
\mbox{in the subgame of level $i$ of }\wGs(r_{i^*_{\zeta_{j^*-2}}}(\zeta_{j^*
-2}),\name{\bbQ}_{\zeta_{j^*-2}})\mbox{ after the two players}\ \\ 
\mbox{played }\langle p^i_{t^i_0}(\zeta_{j^*-2}),q^i_{t^i_{m_{j^*-1}-1}}(
\zeta_{j^*-2})\rangle\mbox{ as the conditions attached to }(t^i_0)_{
\zeta_{j^*-2}}\mbox{ '',}  
  \end{array}\]
and finally $p^i_{t^i_{m_{j^*-1}}}(\xi)=p^i_{t^i_0}(\xi)$ for
$\xi\in\Dom(p^i_{t^i_{m_{j^*-1}}})\setminus (\zeta_{j^*-2}+1)$. Again, a
straightforward induction shows that  
\begin{enumerate}
\item[$(\odot)_{14}^{m_{j^*-1}}$] if $t'=\proj^{(w_i,i+1)}_{(w_{i-1},i)}
\big(t^i_{m_{j^*-1}}\big)$, then $p^i_{t^i_{m_{j^*-1}}}\geq q^{i-1}_{t'}$. 
\end{enumerate}
The condition $p^i_{t^i_{m_{j^*-1}}}$ is played by Generic at stage
$m_{j^*-1}$ and from now till step $m'$ she plays applying procedure 
$\proc_{j^*-1}(m_{j^*-1},m',t^i_{m_{j^*-1}},p^i_{t^i_{m_{j^*-1}}})$.   

Then, if only $m'<m_{j^*-2}$, Generic takes $m''$ such that $[m',m'')$ is
the neighbourhood of $t^i_{m'}$ at $\zeta_{j^*-1}$. She defines
$p^i_{t^i_{m'}}\in\bbP_\gamma$ like $p^i_{t^i_{m_{j^*-1}}}$, so
\[\Dom\big(p^i_{t^i_{m'}}\big)=w_i\cup\Dom(r_{i-1})\cup\big(
\Dom(q^i_{t^i_{m'-1}})\cap\zeta_{j^*-2}\big)\]
and $p^i_{t^i_{m'}}\rest\zeta_{j^*-2}=q^i_{t^i_{m'-1}}\rest\zeta_{j^*-2}$,
and $p^i_{t^i_{m'}}(\zeta_{j^*-2})$ is a $\bbP_{\zeta_{j^*-2}}$--name for an
element of $\name{\bbQ}_{\zeta_{j^*-2}}$ such that   
\[\begin{array}{l}
p^i_{t^i_{m'}}\rest\zeta_{j^*-2}\forces_{\bbP_{\zeta_{j^*-2}}}\mbox{``
$p^i_{t^i_{m'}}(\zeta_{j^*-2})$ is the inning of Generic given by
$\name{\st}_{\zeta_{j^*-2}}$ in }\\ 
\quad \mbox{the subgame of level $i$ of $\wGs(r_{i^*_{\zeta_{j^*-2}}}(
\zeta_{j^*-2}),\name{\bbQ}_{\zeta_{j^*-2}})$ after the two players}\\
\quad\mbox{played }\langle p^i_{t^i_0}(\zeta_{j^*-2}),q^i_{t^i_{m_{j^*-1}
-1}}(\zeta_{j^*-2})\rangle,\langle p^i_{t^i_{m_{j^*-1}}}(\zeta_{j^*-2}),
q^i_{t^i_{m'-1}}(\zeta_{j^*-2})\rangle\mbox{ as}\\
\quad\mbox{the conditions attached to }(t^i_0)_{\zeta_{j^*-2}},
(t^i_{m_{j^*-1}})_{\zeta_{j^*-2}},\mbox{ respectively '',}
 \end{array}\]
and $p^i_{t^i_{m'}}(\xi)=p^i_{t^i_0}(\xi)$ for $\xi\in\Dom(p^i_{t^i_{m'}})
\setminus (\zeta_{j^*-2}+1)$. Like before, 
\begin{enumerate}
\item[$(\odot)_{14}^{m'}$] if $t'=\proj^{(w_i,i+1)}_{(w_{i-1},i)}
\big(t^i_{m'}\big)$, then $p^i_{t^i_{m'}}\geq q^{i-1}_{t'}$. 
\end{enumerate}
The condition $p^i_{t^i_{m'}}$ is played by Generic at stage $m'$ and from
now till step $m''$ she plays applying procedure $\proc_{j^*-1}(m',m'',
t^i_{m'},p^i_{t^i_{m'}})$. Following this pattern untill the players get to 
level $m_{j^*-2}$ (i.e., at the steps from $0$ to $m_{j^*-2}-1$) results in 
defining the procedure $\proc_{j^*-2}(0,m_{j^*-2},t^i_0,p^i_{t^i_0})$.

Suppose we have defined the procedure $\proc_j$, $0<j<j^*$, and we are going
to define $\proc_{j-1}(0,m_{j-1},t^i_0,p^i_{t^i_0})$ following the pattern
presented above. We pick $m'>m_j$ such that $[m_j,m')$ is the neighbourhood
of $t^i_{m_j}$ at $\zeta_j$ and we define $p^i_{t^i_{m_j}}\in\bbP_\gamma$ so
that 
\begin{itemize}
\item $\Dom(p^i_{t^i_{m_j}})=w_i\cup\Dom(r_{i-1})\cup\big(\Dom(q^i_{t^i_{
m_j-1}})\cap\zeta_{j-1}\big)$, 
\item $p^i_{t^i_{m_j}}\rest\zeta_{j-1}=q^i_{t^i_{m_j-1}}\rest\zeta_{j-1}$,
\item $p^i_{t^i_{m_j}}\rest \zeta_{j-1}$ forces in $\bbP_{\zeta_{j-1}}$ that
\begin{quotation}
`` $p^i_{t^i_{m_j}}(\zeta_{j-1})$ is the inning of Generic given to her by
$\name{\st}_{\zeta_{j-1}}$ in the subgame of level $i$ of $\wGs(r_{i^*_{
\zeta_{j-1}}}(\zeta_{j-1}),\name{\bbQ}_{\zeta_{j-1}})$ after the two players
played $\langle p^i_{t^i_0}(\zeta_{j-1}),q^i_{t^i_{m_j-1}}(\zeta_{j-1})
\rangle$ as the conditions attached to $(t^i_0)_{\zeta_{j-1}}$ '',
\end{quotation}
\item $p^i_{t^i_{m_j}}(\xi)=p^i_{t^i_0}(\xi)$ for $\xi\in\Dom(p^i_{
t^i_{m_j}})\setminus (\zeta_{j-1}+1)$.
\end{itemize}
Plainly, $(*)^{m_j}_{14}$ holds and the condition $p^i_{t^i_{m_j}}$ is
played by Generic as her inning associated with $t^i_{m_j}$. Antigeneric
answers with $q^i_{t^i_{m_j}}\geq p^i_{t^i_{m_j}}$, and then Generic follows
procedure $\proc_j(m_j,m',t^i_{m_j},p^i_{t^i_{m_j}})$ to determine her
innings for $m\in (m_j,m')$. If $m'<m_{j-1}$, Generic picks $m''>m'$ such
that $[m',m'')$ is the neighbourhood of $t^i_{m'}$ at $\zeta_j$ and she
defines $p^i_{t^i_{m'}}\in\bbP_\gamma$ like before, and then she follows the
procedure $\proc_j(m',m'',t^i_{m'},p^i_{t^i_{m'}})$, and so on. After
arriving to $m_{j-1}-1$ Generic defined the procedure $\proc_{j-1}(0,
m_{j-1},t^i_0,p^i_{t^i_0})$. 

Finally, the procedure $\proc_0(0,m_0,t^i_0,p^i_{t^i_0})$ describes the
instructions given to Generic by our strategy $\st$ in the subgame of level
$i$ of $\wGs(p,\bbP_\gamma)$. Now, for $\xi\in w_i$, Generic defines
$\name{\bar{p}}_{i,\xi}=\langle \name{p}^\eta_{i,\xi}:\eta\in
\max(s_{i,\xi})\rangle$ and $\name{\bar{q}}_{i,\xi}=\langle\name{q}^\eta_{i,
\xi}:\eta\in\max(s_{i,\xi})\rangle$ so that $\name{p}^\eta_{i,\xi},
\name{q}^\eta_{i,\xi}$ are $\bbP_\xi$--names for elements of
$\name{\bbQ}_\xi$ and the demands in $(\odot)_6+(\odot)_{11}$ are
satisfied. It should be clear that the choice of $\name{\bar{p}}_{i,\xi},
\name{\bar{q}}_{i,\xi}$ is possible --- note that by niceness of
$\name{\st}_\xi$ (and the last requirement in \ref{nice}(1)) we easily
justify that $(\odot)_{10}$ holds. To finish stage $i$, Generic picks
$r_i\in\bbP_\gamma$ essentially by conditions $(\odot)_2+(\odot)_{12}$. 
\medskip

We are going to argue that $\st$ is a winning strategy for Generic in
$\tGs(p,\bbP_\gamma)$. To this end suppose that $\langle (\cT_i,\bar{p}^i, 
\bar{q}^i):i<\omega\rangle$ is the result of a play of
$\tGs(p,\bbP_\gamma)$ in which Generic used $\st$, and the side objects
constructed at each stage $i<\omega$ are
\begin{enumerate}
\item[$(\boxdot)$] $w_i,\cT_i,r_i,\bar{p}^i,\bar{q}^i,k_i,K_\xi,i^*_\xi,
\name{\st}_\xi,s_{i,\xi},\name{\bar{p}}_{i,\xi},\name{\bar{q}}_{i,\xi}$ for
$\xi\in w_i$,  
\end{enumerate}
and they satisfy the requirements $(\odot)_0$---$(\odot)_{12}$. 

We define a condition $q\in\bbP_\gamma$ as follows. Let $\Dom(q)=
\bigcup\limits_{i<\omega} w_i=\bigcup\limits_{i<\omega}\Dom(r_i)$ and for
$\xi\in\Dom(q)$ let $q(\xi)$ be a $\bbP_\xi$--name for a condition in
$\name{\bbQ}_\xi$ such that
\begin{enumerate}
\item[$(\odot)^\xi_{15}$] $\forces_{\bbP_\xi}\mbox{`` }q(\xi)\geq
r_{i^*_\xi}(\xi)\mbox{ and } q(\xi)\forces_{\name{\bbQ}_\xi} \big(\forall
i\geq i^*_\xi\big)\big(\exists \nu\in \max(s_{i,\xi})\big)\big(
\name{q}^\nu_{i,\xi}\in\Gamma_{\name{\bbQ}_\xi}\big)$ ''.
\end{enumerate}
Clearly $q$ is well defined (remember $(\odot)_6$) and $q\geq p$ (remember 
$(\odot)_1+(\odot)_2$). Also $q\geq r_i$ for all $i<\omega$. 

We will show that for each $i<\omega$ the family $\{q^i_t:t\in T_i\ \&\
\rk_i(t)=\gamma\}$ is predense above $q$ (and this clearly implies that
Generic won the play). So suppose $q^+\geq q$, $i<\omega$ and $w_i\cup\{ 
\gamma\}=\{\xi_0,\xi_1,\ldots,\xi_i,\xi_{i+1}\}$ (the increasing
enumeration, so $\xi_0=0$, $\xi_{i+1}=\gamma$). By induction on $j\leq i$ we
choose an increasing sequence $\langle q^+_j:j\leq i\rangle\subseteq
\bbP_\gamma$ and we also define $t\rest\xi_j+1$. First, by the choice of
$q(0)$ there is $\nu\in\max(s_{0,\xi_0})$ such that the conditions $q^+(0)$
and $q^\nu_{i,\xi_0}=\name{q}^\nu_{i,\xi_0}$ are compatible. Let $(t)_0=\nu$ 
(this defines $t\rest\xi_1$). Let $q^+_0\in \bbP_{\xi_1}$ be such that
$q^+_0(0)$ is stronger than both $q^+(0)$ and $q^\nu_{i,\xi_0}$, and
$q^+_0\rest (\xi_0,\xi_1)=q^+\rest (\xi_0,\xi_1)$. It follows from 
$(\odot)_{10}+(\odot)_{12}$ that 
\begin{enumerate}
\item[$(\odot)^0_{16}$] $\big(\forall t'\in T_i\big)\big([\rk_i(t')=\gamma\
\&\ t'\rest\xi_1=t\rest\xi_1]\ \Rightarrow\ q^i_{t'}\rest\xi_1\leq
q^+_0\big)$. 
\end{enumerate}
Suppose that $j<i$ and we have defined $t\rest\xi_{j+1}\in T_i$ and a
condition $q^+_j\in \bbP_{\xi_{j+1}}$ such that $q^+_j\geq
q^+\rest\xi_{j+1}$ and 
\begin{enumerate}
\item[$(\odot)^j_{16}$] $\big(\forall t'\in T_i\big)\big([\rk_i(t')=\gamma\
\&\ t'\rest\xi_{j+1}=t\rest\xi_{j+1}]\ \Rightarrow\ q^i_{t'}\rest\xi_{j+1}
\leq q^+_j\big)$. 
\end{enumerate}
Necessarily, by $(\odot)_{15}^{\xi_{j+1}}$, 
\[q^+_j\forces_{\bbP_{\xi_{j+1}}}\mbox{`` }\big(\exists\nu\in
\max(s_{i,\xi_{j+1}})\big)\big(\name{q}^\nu_{i,\xi_{j+1}}, q^+(\xi_{j+1})
\mbox{ are compatible }\big)\mbox{ ''}\]
so we may choose $\nu\in\max(s_{i,\xi_{j+1}})$ and a condition $q_{j+1}\in
\bbP_{\xi_{j+1}}$ stronger than $q^+_j$ such that  
\[q_{j+1}\forces_{\bbP_{\xi_{j+1}}}\mbox{`` }\name{q}^\nu_{i,\xi_{j+1}},
q^+(\xi_{j+1})\mbox{ are compatible ''}.\] 
Let $(t)_{\xi_{j+1}}=\langle *\ldots*\rangle\conc\nu$ (thus $t\rest(\xi_{j+
1}+1)$ has been defined) and let $q^+_{j+1}\in \bbP_{\xi_{j+2}}$ be such
that $q^+_{j+1}\rest\xi_{j+1}=q_{j+1}$, 
\begin{enumerate}
\item[$(\odot)_{17}$] $q^+_{j+1}\rest\xi_{j+1}\forces_{\bbP_{\xi_{j+1}}}$``
$q^+_{j+1}(\xi_{j+1})\geq \name{q}^\nu_{i,\xi_{j+1}}\quad \&\quad q^+_{j+1}
(\xi_{j+1})\geq q^+(\xi_{j+1})$ '', 
\end{enumerate}
and $q^+_{j+1}\rest(\xi_{j+1},\xi_{j+2})=q^+\rest(\xi_{j+1},\xi_{j+2})$.
It follows from $(\odot)_{11}+(\odot)_{16}$ that 
\begin{enumerate}
\item[$(\odot)_{18}$] if $t'\in T_i$, $\rk_i(t')=\gamma$, $t'\rest\xi_{j
+2}=t\rest\xi_{j+2}$, then $q^+_{j+1}\rest \xi_{j+1}
\forces_{\bbP_{\xi_{j+1}}} q^i_{t'}(\xi_{j+1})\leq \name{q}^\nu_{i,
\xi_{j+1}}$. 
\end{enumerate}
We may use $(\odot)_{17}+(\odot)_{18}$ and then $(\odot)_{12}$ to argue that
$(\odot)_{16}^{j+1}$ holds true.

Finally look at $t=t\rest\xi_{i+1}$ and $q^+_{i+1}$. 
\end{proof}

\begin{definition}
\label{barst}
Suppose that $\bbP$ is a forcing notion with $\odot_n$--property and
$\bar{\st}=\langle\st_p:p\in\bbP\rangle$, where $\st_p$ is a winning
strategy for Generic in $\wGs(p,\bbP)$. Such $\bar{\st}$ will be called {\em
an $\odot_n$--strategy system for $\bbP$}.
\begin{enumerate}
\item We say that a finite set $Q$ of conditions is {\em an
$\bar{\st}$--front above $p$\/} provided that there is a partial play
$\langle s_j,\bar{p}^j,\bar{q}^j:j\leq i\rangle$ of $\wGs(p,\bbP)$ in
which Generic uses $\st_p$ and
\begin{itemize}
\item if $\max(s_i)=\langle\eta^i_k:k<K\rangle$ is the enumeration played by
Antigeneric at stage $i$ after Generic put $s_i$, then $Q=
\{q^i_{\eta^i_k}:k<K\}$.    
\end{itemize}
\item For a condition $p\in\bbP$ we define a game $\jeszcze(p,\bbP)$ as
follows. A play of $\jeszcze(p,\bbP)$ lasts $\omega$ moves and in the course
of the play a sequence 
\begin{enumerate}
\item[$(\boxtimes)$] \qquad $\big\langle s_i,\bar{\eta}^i,\bar{p}^i,
\bar{Q}^i,\bar{q}^i:i<\omega\big\rangle$  
\end{enumerate}
is constructed. At a stage $i<\omega$ of the play, 
\begin{itemize}
\item first Generic chooses a finite $n$--ary tree $s_i$ such that the
demand $(\alpha)$ of \ref{da2}(1) holds, and then 
\item Antigeneric picks an enumeration $\bar{\eta}^i=\langle\eta^i_\ell:
\ell<k_i\rangle$ of $\max(s_i)$.  
\end{itemize}
Now the two players start a subgame of length $k_i$ and they choose
successive terms of a sequence $\langle p^i_{ \eta^i_\ell}, Q^i_{
\eta^i_\ell}:\ell<k_i\rangle$. At a stage $\ell<k_i$ of the subgame, 
\begin{itemize}
\item first Generic picks a condition $p^i_{\eta^i_\ell}\in\bbP$ such that 
\smallskip

if $j<i$, $\nu\in \max(s_j)$ and $\nu\vtl\eta^i_\ell$, then
$q^j_\nu\leq p^i_{\eta^i_\ell}$ and $p\leq p^i_{\eta^i_\ell}$,  

\item and then Antigeneric picks an $\bar{\st}$--front $Q^i_{\eta^i_\ell}$
above $p^i_{\eta^i_\ell}$.  
\end{itemize}
After the subgame is completed, 
\begin{itemize}
\item Antigeneric chooses $\bar{q}^i=\langle q^i_\eta:\eta\in\max(s_i)
\rangle$ so that $q^i_\eta\in Q^i_\eta$ for $\eta\in\max(s_i)$. 
\end{itemize}
Finally, Generic wins a play $(\boxtimes)$ if and only if  
\begin{enumerate}
\item[$(\circledast)$] there is a condition $q\geq p$ such that for every
$i<\omega$ the family $\{q^i_\eta:\eta\in\max(s_i)\}$ is predense above $q$.  
\end{enumerate}
\item Similarly to \ref{nice}(1) we define when a strategy $\st$ of Generic
in $\jeszcze(p,\bbP)$ is {\em nice for an infinite set $K\subseteq\omega$}. 
\item We say that the forcing notion $\bbP$ has {\em the uniformly nice
$(\odot)^{\bar{\st}}_n$--property\/} if for every $p\in\bbP$ and an infinite
set $K\subseteq\omega$ Generic has a nice for $K$ winning strategy in
$\jeszcze(p,\bbP)$.  
\end{enumerate}
\end{definition}

\begin{observation}
\begin{enumerate}
\item If $\bbP$ has the $\oplus_n$--property, then it  has the uniformly
nice $(\odot)^{\bar{\st}}_n$--property for some $\odot_n$--strategy system 
$\bar{\st}$. 
\item The uniform $n$--Sacks forcing notion $\bbQ_n$ has the uniformly nice 
$(\odot)^{\bar{\st}}_n$--property for some $\odot_n$--strategy system
$\bar{\st}$.  
\end{enumerate}
\end{observation}

We do not know if the $n$--Silver forcing is equivalent to a forcing with
the uniformly nice $(\odot)^{\bar{\st}}_n$--property (for some $\bar{\st}$).   

\begin{theorem}
\label{getnloc}
Assume that $\bar{\bbQ}=\langle\bbP_\xi,\name{\bbQ}_\xi:\xi<\gamma\rangle$
is a CS iteration and $\name{\bar{\st}}^\xi$ are $\bbP_\xi$--names such that
for every $\xi<\gamma$,  
\[\begin{array}{ll}
\forces_{\bbP_\xi}&\mbox{`` $\name{\bbQ}_\xi$ has the $\odot_n$--property
and }\name{\bar{\st}}^\xi\mbox{ is a $\odot_n$--strategy system for
$\name{\bbQ}_\xi$ and}\\
&\ \ \name{\bbQ}_\xi\mbox{ has the uniformly nice
$(\odot)_n^{\name{\bar{\st}}^\xi}$--property ''.}
  \end{array}\]  
Then $\bbP_{\gamma}=\lim(\bar{\bbQ})$ has the $n$--localization property.
\end{theorem}

\begin{proof}
The following combinatorial observation can be shown by an easy induction. 

\begin{claim}
\label{cl1}
Let $M<\omega$. Suppose that for each $m<M$ we are given $k_m<\omega$ and a
set $A_m\subseteq {}^{k_m}\omega$ of size $M!$. Then there is a sequence
$\langle\sigma_m:m<M \rangle\in \prod\limits_{m<M} A_m$ such that 
\[\big(\forall m<m'<M\big)\big(\sigma_m,\sigma_{m'}\mbox{ are incompatible
}\big).\] 
\end{claim}

Let $\name{\tau}$ be a $\bbP_\gamma$--name for an element of
$\baire$. Without loss of generality we may assume that 
\begin{enumerate}
\item[$(\otimes)_0$]\qquad $\forces_{\bbP_\gamma}(\forall\alpha<\gamma)
  (\name{\tau}\notin \bV^{\bbP_\alpha})$. 
\end{enumerate}

\begin{claim}
\label{cl2}
Let $K,M<\omega$, $\xi<\gamma$, $p\in\bbP_\gamma$. Then there are $N>M$,
$q^*\in\bbP_\xi$ and a sequence $\langle\sigma_\ell,q_\ell:\ell<L\rangle$
such that 
\begin{enumerate}
\item[(a)] $\sigma_\ell\in {}^N\omega$ and $|\{\sigma_\ell\rest [M,N):
\ell<L\}|>K$, 
\item[(b)] $q_\ell\in\bbP_\gamma$, $p\leq q_\ell$, $q_\ell\rest\xi=q^*$ and
$q_\ell\forces_{\bbP_\gamma}\name{\tau}\rest N=\sigma_\ell$,
\item[(c)] $q^*\forces_{\bbP_\xi}$`` $\{q_\ell(\xi):\ell<k\}$ is an
  $\name{\bar{\st}}^\xi$--front above $p(\xi)$ ''. 
\end{enumerate}
\end{claim}

\begin{proof}[Proof of the Claim]
Let $\name{\bbP}_{\xi\gamma}$ be defined as at the beginning of the proof of
Lemma \ref{treelem}. Suppose that $G_\xi\subseteq\bbP_\xi$ is generic over
$\bV$, $p\rest\xi\in G_\xi$, and let us work in $\bV[G_\xi]$ for a while. 
Then $\name{\bbP}_{\xi\gamma}[G_\xi]$ is a dense subset of the limit
$\bbP^{\xi\gamma}$ of a CS iteration of forcing notions with
$\odot_n$--property, so we may use the proof of \ref{seciter}. Consider a
play $\langle(\cT_i,\bar{p}^i, \bar{q}^i):i<\omega\rangle$ of $\tGs(p
\rest[\xi,\gamma),\bbP^{\xi\gamma})$ in which 
\begin{itemize}
\item Generic follows exactly the strategy described in the proof of
\ref{seciter}, where on the coordinate $\xi$ the strategy
$\name{\st}^\xi_{p(\xi)}[G_\xi]$ is used and  
\item each condition $q^i_t$ played by Antigeneric (for $t\in T_i$,
$\rk_i(t)=\gamma$) is from $\name{\bbP}_{\xi\gamma}[G_\xi]$ and decides the
value of $\name{\tau}\rest (M+i)$, say $q^i_t\forces \name{\tau}\rest
(M+i)=\sigma^i_t$. 
\end{itemize}
Let $q\in\bbP^{\xi,\gamma}$ be the condition defined by $(\odot)^\xi_{15}$
(at the end of the proof of \ref{seciter}), so it witnesses that Generic won
the play, and $\{q^i_t:t\in T_i\ \&\ \rk_i(t)=\gamma\}$ is predense above
$q$ (for each $i<\omega$). It follows from $(\otimes)_0$ that for some
$i<\omega$ we have 
\[K<\big|\big\{\sigma^i_t\rest [M,M+i):t\in T_i\ \&\ \rk_i(t)=\gamma\big\}
  \big|.\]
Also it follows from the description of Generic's strategy in $\tGs(p\rest[
\xi,\gamma),\bbP^{\xi\gamma})$ that the family $\{q^i_t(\xi):
t\in T_i\ \&\ \rk_i(t)=\gamma\}$ is an $\name{\bar{\st}}^\xi[G_\xi]$--front
above $p(\xi)$. 

Now, $T_i,\langle\sigma^i_t,q^i_t:t\in T_i\ \&\ \rk_i(t)=\gamma\rangle\in
\bV$, so we may pick a condition $q^*\in\bbP_\xi$ stronger than $p\rest\xi$
which forces that these objects have the properties described above. Let
$\langle t_\ell:\ell<L\rangle$ be an enumeration of $\{t\in T_i:\rk_i(t)=
\gamma\}$ and $\sigma_\ell=\sigma^i_{t_\ell}$, $q_\ell=q^*\conc
q^i_{t^i_\ell}$ (for $\ell<L$).     
\end{proof}

Let $p\in\bbP_\gamma$. Following the procedure described in the proof of
Theorem \ref{seciter} construct a sequence 
\[\big\langle w_i,\cT_i,\bar{p}^i,\bar{q}^i,\bar{q}^{i,*},\bar{\sigma}^i,
r_i,k_i,\alpha_i,M_i,\bar{N}^i, i^*_\xi, \name{\st}_\xi,s_{i,\xi},
\bar{\eta}_{i,\xi},\name{\bar{p}}_{i,\xi},\name{\bar{Q}}_{i,\xi},
\name{\bar{q}}_{i,\xi}:\ i\!<\!\omega,\ \xi\!\in\! w_i\big\rangle\]   
such that the following conditions are satisfied. 
\begin{enumerate}
\item[$(\otimes)_1$] The demands formulated in $(\odot)_0$--$(\odot)_3$ and 
$(\odot)_7$--$(\odot)_{12}$ of the proof of \ref{seciter}. 
\item[$(\otimes)_2$] For $\xi\in\bigcup\limits_{i<\omega} w_i$ we have
$i^*_\xi=\min(\{i:\xi\in w_i\})\leq\min(K_\xi)$, and $\name{\st}_\xi$ is a 
$\bbP_\xi$--name for a winning strategy of Generic in $\jeksi(r_{i^*_\xi}
(\xi),\name{\bbQ}_\xi)$ which is nice for $\{k\in\omega:k+i^*_\xi\in K_\xi
\}$.      
\item[$(\otimes)_3$] If $\xi\in w_i$, then $s_{i,\xi}\subseteq
\bigcup\limits_{j\leq i+1-i^*_\xi} {}^j(n+1)$ is an $n$--tree,
$\bar{\eta}_{i,\xi}$ is the enumeration of $\max(s_{i,\xi})$ in the
$<^*_\chi$--increasing order, and $\name{\bar{p}}_{i,\xi}=\langle
\name{p}^\eta_{i,\xi}:\eta\in\max(s_{i,\xi})\rangle$,
$\name{\bar{Q}}_{i,\xi}=\langle\name{Q}^\eta_{i,\xi}:\eta\in
\max(s_{i,\xi})\rangle$, $\name{\bar{q}}_{i,\xi}=\langle\name{q}^\eta_{i, 
\xi}:\eta\in\max(s_{i,\xi})\rangle$ are $\bbP_\xi$--names for systems
indexed by $\max(s_{i,\xi})$.     
\item[$(\otimes)_4$] For each  $\xi\in \bigcup\limits_{i<\omega} w_i$,  
\[\begin{array}{r}
\forces_{\bbP_\xi}\mbox{`` }\langle s_{i,\xi},\bar{\eta}_{i,\xi},
\name{\bar{p}}_{i,\xi},\name{\bar{Q}}_{i,\xi},\name{\bar{q}}_{i,\xi}:
i^*_\xi\leq i<\omega\rangle\mbox{ is a legal play of }\jeksi(r_{i^*_\xi}
(\xi),\name{\bbQ}_\xi)\\ 
\mbox{ in which Generic uses $\name{\st}_\xi$ ''.}
\end{array}\]
\item[$(\otimes)_5$] $\bar{N}^i=\langle N^i_t:t\in T_i\ \&\ \rk_i(t)=\gamma
\rangle$, $M_i<N^i_t<M_{i+1}<\omega$, $M_0=0$, and $\bar{\sigma}^i=\langle
\sigma^i_{t,\ell}:t\in T_i\ \&\ \rk_i(t)=\gamma\ \&\ \ell<k_i!\rangle$,
$\sigma^i_{t,\ell}\in {}^{N^i_t}\omega$, and if $\ell<\ell'<k_i!$, $t\in
T_i$, $\rk_i(t)=\gamma$, then $\sigma^i_{t,\ell}\rest [M_i, N^i_t)\neq
\sigma^i_{t,\ell'}\rest [M_i, N^i_t)$.

\item[$(\otimes)_6$] $\alpha_i=\max(w_i)$, $\bar{q}^{i,*}=\langle
q^{i,*}_{t,\ell}:t\in T_i\ \&\ \rk_i(t)=\gamma\ \&\ \ell<k_i!\rangle$, and
for $t\in T_i$, $\rk_i(t)=\gamma$, $\ell<k_i!$ we have $q^{i,*}_{t,\ell}
\rest\alpha_i=q^i_t\rest\alpha_i$, $q^{i,*}_{t,\ell}\rest\alpha_i
\forces_{\bbP_{\alpha_i}} q^{i,*}_{t,\ell}(\alpha_i)\in
\name{Q}^\eta_{i,\alpha_i}$, where $\eta\in\max(s_{i,\alpha_i})$ is the end 
segment of $(t)_{\alpha_i}$ of length $i+1-i^*_\xi$, and $q^{i,*}_{t,\ell}
\forces_{\bbP_\gamma}\name{\tau}\rest N^i_t=\sigma^i_{t,\ell}$. 
\item[$(\otimes)_7$] If $t_0,t_1\in T_i$, $\rk_i(t_0)=\rk_i(t_1)=\gamma$,
$t_0\neq t_1$, then for some $\ell_0,\ell_1<k_i!$ we have that $q^i_{t_0}
=q^{i,*}_{t_0,\ell_0}$, $q^i_{t_1}=q^{i,*}_{t_1,\ell_1}$  and  the
sequences $\sigma^i_{t_0,\ell_0}\rest [M_i,N^i_{t_0})$ and $\sigma^i_{t_1,
\ell_1}\rest [M_i,N^i_{t_1})$  are incomparable 
\end{enumerate}
To guarantee demands $(\otimes)_1$--$(\otimes)_4$ we follow exactly the
lines of the proof of \ref{seciter}, to get $(\otimes)_5+(\otimes)_6$ we
use Claim \ref{cl2} and we ensure $(\otimes)_7$ by Claim \ref{cl1}.  

After the construction is carried out define a condition $q\in\bbP_\gamma$
in a manner similar to that in the proof or \ref{seciter}: $\Dom(q)=
\bigcup\limits_{i<\omega} w_i=\bigcup\limits_{i<\omega}\Dom(r_i)$ and for
$\xi\in\Dom(q)$ we let $q(\xi)$ be a $\bbP_\xi$--name for a condition in
$\name{\bbQ}_\xi$ such that
\begin{enumerate}
\item[$(\otimes)^\xi_8$] $\forces_{\bbP_\xi}\mbox{`` }q(\xi)\geq
r_{i^*_\xi}(\xi)\mbox{ and } q(\xi)\forces_{\name{\bbQ}_\xi} \big(\forall
i\geq i^*_\xi\big)\big(\exists \nu\in \max(s_{i,\xi})\big)\big(
\name{q}^\nu_{i,\xi}\in\Gamma_{\name{\bbQ}_\xi}\big)$ ''.
\end{enumerate}
As in \ref{seciter} one argues that 
\begin{enumerate}
\item[$(\otimes)_9$] for each $i<\omega$ the family $\{q^i_t:t\in T_i\ \&\
\rk_i(t)=\gamma\}$ is predense above $q$.
\end{enumerate}
Now we choose a tree $T\subseteq {}^{\omega{>}}\omega$ such that 
\[(\forall f\in [T])(\forall i<\omega)(\exists t\in T_i)(\exists \ell<k_i)(
q^i_t=q^{i,*}_{t,\ell}\ \&\ \sigma^i_{t,\ell}\vtl f).\] 
Plainly, $T$ is an $n$--ary tree and $q\forces_{\bbP_\gamma}\name{\tau}\in
[T]$. 
\end{proof}

\begin{remark}
\label{finrem}
After analyzing the proof of Theorem \ref{getnloc} one may notice that the
following can be shown by the same proof. 
\begin{enumerate}
\item[] {\em Assume that $\gamma$ is a limit ordinal and $\bar{\bbQ}=\langle
\bbP_\xi,\name{\bbQ}_\xi:\xi< \gamma\rangle$ is a CS iteration such that for
every $\xi<\gamma$
\begin{itemize}
\item $\forces_{\bbP_\xi}$`` $\name{\bbQ}_\xi$ has the nice
  $\odot_n$--property '', 
\item $\bbP_\xi$ has the $n$--localization property.
\end{itemize}
Then $\bbP_{\gamma}=\lim(\bar{\bbQ})$ has the $n$--localization property.}
\end{enumerate}
(The assumption that $\gamma$ is limit allows us to make sure in the
construction that $i<\min(K_{\alpha_i})$ for all $i<\omega$.)
\end{remark}

\begin{problem}
\begin{enumerate}
\item Can the implications in Observation \ref{easyobs} be reversed? 
What if we restrict ourselves to (s)nep forcing notions or even Suslin$^+$ ?  
\item Assume that $\bbP$ has the $\odot_n$--property. Is it equivalent to a
  forcing notion with the uniformly nice $(\odot)_n^{\bar{\st}}$--property
  (for some $\odot_n$--strategy system $\bar{\st}$)~? Again, we may allow
  restrictions to nice forcing notions. 
\end{enumerate}
\end{problem}


\end{document}